\documentclass[12pt]{amsart}
\usepackage{amsmath}
\usepackage{amssymb}
\usepackage{latexsym}

\newtheorem{prop}{Proposition}[section]
\newtheorem{thm}[prop]{Theorem}

\newtheorem{cor}[prop]{Corollary}
\newtheorem{conj}{Conjecture}

\theoremstyle{definition}
\newtheorem{df}{Definition}
\newtheorem*{remk}{Remark} 

\newcommand{\Zint}{\mathbb {Z}}    
     
\newcommand{\Rea}{\mathbb {R}}      
\newcommand{\Cplx}{\mathbb {C}}     

\newcommand{\halmos}{\rule{5pt}{5pt}}

\begin{document}

\title[Heun equation {\rm III}]
{The Heun equation and the Calogero-Moser-Sutherland system {\rm III}: the finite-gap property and the monodromy}
\author{Kouichi Takemura}
\address{Department of Mathematical Sciences, Yokohama City University, 22-2 Seto, Kanazawa-ku, Yokohama 236-0027, Japan.}
\email{takemura@yokohama-cu.ac.jp}

\begin{abstract}
\noindent
A new approach to the finite-gap property for the Heun equation is constructed. The relationship between the finite-dimensional invariant space and the spectral curve is clarified. 
The monodromies are calculated and are expressed as hyperelliptic integrals.
Applications to the spectral problem for the $BC_1$ Inozemtsev model are obtained.
\end{abstract}

\maketitle

%
\section{Introduction}

The $BC_1$ Inozemtsev model is a one-particle quantum mechanics model whose Hamiltonian is given as
\begin{equation}
H= -\frac{d^2}{dx^2} + \sum_{i=0}^3 l_i(l_i+1)\wp (x+\omega_i),
\label{Ino}
\end{equation}
where $\wp (x)$ is the Weierstrass $\wp$-function with periods $(1,\tau )$, $\omega _0=0$, $\omega_1=1/2$, $\omega_2=-(\tau+1)/2$ and $\omega_3=\tau/2$ are half-periods, and $l_i$ $(i=0,1,2,3)$ are coupling constants.
 This model is a one-particle version of the $BC_N$ Inozemtsev system \cite{Ino}, which is known to be the universal quantum integrable system with $B_N$ symmetry \cite{Ino,OOS}.
The $BC_N$ Calogero-Moser-Sutherland systems are special cases of $BC_N$ Inozemtsev systems.

It is known that there is a correspondence between the $BC_1$ Inozemtsev model and the Heun equation. Here the Heun equation is a standard canonical form of the Fuchsian equation with four singularities defined in (\ref{Heun}) (see section \ref{sec:Heun}).

The potential $ \sum_{i=0}^3 l_i(l_i+1)\wp (x+\omega_i)$ has also arisen from a different source, the theory of elliptic solitons. Around 1987-92 Treibich and Verdier \cite{TV} found and showed that $ \sum_{i=0}^3 l_i(l_i+1)\wp (x+\omega_i)$ is an algebro-geometric finite-gap potential iff $l_i \in \Zint $ for all $i \in \{ 0,1,2,3 \}$. For this reason the function $ \sum_{i=0}^3 l_i(l_i+1)\wp (x+\omega_i)$ is called the Treibich-Verdier potential. 
Subsequently several others have produced results on this subject, e.g., Gesztesy and Weikard \cite{GW,Wei}, Smirnov \cite{Smi} and so forth. (For more details see the references in \cite{GW,Wei,Smi}.)

In this paper a new approach to the Treibich-Verdier potential is presented, along with applications to the spectral problem for the $BC_1$ Inozemtsev model.
In particular the relationship between some finite-dimensional space of elliptic functions and the spectral curve is clarified and a new determinantal formula for a commuting operator of odd degree is obtained. 
Here the existence of a commuting operator of odd degree is equivalent to the algebro-geometric finite-gap property.

The reader may find the author's background in this subject to be helpful in elucidating the motivations for this research.
In \cite{Tak1} I proposed and developed the Bethe ansatz method for the Heun equation and applied it to the analysis of the eigenstates and the eigenvalues for the $BC_1$ Inozemtsev model.
A few months after the work \cite{Tak1} I noticed that the function $ \sum_{i=0}^3 l_i(l_i+1)\wp (x+\omega_i)$ is an algebro-geometric finite-gap potential while still being unaware of the results of Treibich and Verdier. An argument based on the Bethe ansatz led me to conjecture that the global monodromies are expressible in terms of hyperelliptic integrals and that they would be applicable for investigating the singularities of the eigenvalues in $\tau$.
Eventually the results of Treibich and Verdier \cite{TV}, Hochstadt \cite{Hoc}, McKean and van Moerbeke \cite{MvM} and Smirnov \cite{Smi} came to my attention.
In particular some formulas for the monodromy are given in \cite{Hoc,MvM} and it seemed natural that they would be applicable for the Treibich-Verdier potential.

We next present a brief sketch of the contents and results of this paper.

Firstly a method of determining whether the function $ \sum_{i=0}^3 l_i(l_i+1)\wp (x+\omega_i)$ $(l_0, l_1, l_2, l_3  \in \Zint )$ is an algebro-geometric finite-gap potential is presented in sections \ref{sec2} and \ref{sec3}. The strategy is as follows: First the product of some two solutions to $H f(x)=Ef(x)$ (with $H$ the Hamiltonian, $E$ the eigenvalue) is shown to be doubly periodic for all $E$, a result which was already obtained in \cite{Tak1} in order to show the completeness of the Bethe ansatz.
From the solution to the equation for that product a differential operator $A$ of odd degree that commutes with the Hamiltonian $H$ is constructed in section \ref{sec:constr}, whence follows the algebro-geometric finite-gap property. The proof is elementary and does not require the advanced techniques of Treibich and Verdier that are based on the theory of Jacobi varieties and the notion of a tangential cover.

In section \ref{sec:invsp} finite-dimensional spaces of elliptic functions are introduced.
Note that these finite-dimensional spaces have already appeared in \cite{Tak1} and they are related to the concept of quasi-exact solvability \cite{Tur,GKO}.
These spaces play an important role in the presentation and investigation of the spectral curve, i.e. the curve determined by the equality between the Hamiltonian $H$ and the commuting operator $A$. As an application we obtain a determinantal formula for the commuting operator in section \ref{sec:det}.
In section \ref{sec:monod} formulas for the monodromy of the solutions to the Inozemtsev model are obtained. These formulas are expressed in the form of hyperelliptic integrals.
In section \ref{sec:exam} examples are used to illustrate the results and the formulas obtained in sections \ref{sec2} and \ref{sec3}.

From the viewpoint of physics it is often desirable to obtain the square-integrable eigenfunctions of the Hamiltonian and the corresponding eigenvalues. In sections \ref{sec2} and \ref{sec3} the eigenfunctions of $H$ without this physical boundary condition are constructed, whereas in section \ref{sec4} the physical boundary condition is included.

It is shown in \cite{Tak2} that an eigenvalue satisfying the physical boundary condition is holomorphic in $\tau$ if $-1<p=\exp (\pi \sqrt{-1} \tau )<1$. In this paper a sufficient condition is provided for the holomorphy of the physical eigenvalues in $\tau \in \Cplx$. This condition involves the convergenec radius of the eigenvalue as a power series in $p$, which can be calculated by the method of perturbation \cite{Tak2}. 

In section \ref{sec:str} the strategy is outlined. Section \ref{sec:Hilbsp} provides the appropriate setting for consideration of the spectral problem in the case $l_0 , l_1 \in \Zint_{\geq 1}$. In particular the Hilbert space ${\bold H}$ and its dense subspace ${\bold D}$ are introduced. In section \ref{sec:condsq} we obtain a sufficient condition for the holomorphy of the physical eigenvalue in $\tau \in \Cplx$ for this case. Here the formulas for the monodromy from section \ref{sec:monod} are used. An example is also presented to illustrate the results.
In section \ref{sec:l10} the case $l_0 \in \Zint_{\geq 1}$ and $l_1=0$ is discussed and in section \ref{sec:l0l10} the case $l_0 = l_1=0$ is treated.

In section \ref{sec:Heun} the relationship between the Heun equation and the Inozemtzev model is reviewed and a translation between the two formalisms is detailed.
Formulas for the monodromy along the cycle around two regular singular points are derived. They are expressed in use of hyperelliptic integrals.

In the appendix definitions and formulas for the elliptic functions are provided.

Throughout this paper it is assumed that $l_i \in \Zint $ $(i=0,1,2,3)$ and $(l_0,l_1,l_2,l_3) \neq (0,0,0,0)$.

\section{Results from \cite{Tak1} and related topics} \label{sec2}

\subsection{Product of the solutions}$ $

Fix the eigenvalue $E$ of the Hamiltonian $H$ (\ref{Ino}) and consider the second-order differential equation
\begin{equation}
(H-E) f(x)= \left( -\frac{d^2}{dx^2} + \sum_{i=0}^3 l_i(l_i+1)\wp (x+\omega_i)-E\right) f(x)=0.
\label{InoEF}
\end{equation}

Let $h(x)$ be  the product of any pair of the solutions to (\ref{InoEF}). Then the function $h(x)$ satisfies the following third-order differential equation:
\begin{align}
& \left( \frac{d^3}{dx^3}-4\left( \sum_{i=0}^3 l_i(l_i+1)\wp (x+\omega_i)-E\right)\frac{d}{dx}-2\left(\sum_{i=0}^3 l_i(l_i+1)\wp '(x+\omega_i)\right) \right) h (x)=0.
\label{prodDE} 
\end{align}

It is known that equation (\ref{prodDE}) has a nonzero doubly periodic solution for all $E$ if $l_i \in \Zint _{\geq 0}$ $(i=0,1,2,3)$. That was first observed by Darboux \cite{Dar}. 
\begin{prop} \cite[Proposition 3.5]{Tak1} \label{prop:prod}
If $l_0, l_1, l_2, l_3 \in \Zint _{\geq 0}$, then equation (\ref{prodDE}) has a nonzero doubly periodic solution $\Xi (x,E)$, which has the expansion
\begin{equation}
\Xi (x,E)=c_0(E)+\sum_{i=0}^3 \sum_{j=0}^{l_i-1} b^{(i)}_j (E)\wp (x+\omega_i)^{l_i-j},
\label{Fx}
\end{equation}
where the coefficients $c_0(E)$ and $b^{(i)}_j(E)$ are polynomials in $E$, they do not have common divisors and the polynomial $c_0(E)$ is monic.
We set $g=\deg_E c_0(E)$. Then the coefficients satisfy $\deg _E b^{(i)}_j(E)<g$ for all $i$ and $j$.
\end{prop}

We can derive the integral formula for the solution $\Lambda(x,E)$ to the equation (\ref{InoEF}) in terms of the doubly periodic function $\Xi(x,E)$, which is obtaind in \cite{Tak1}. Set
\begin{align}
 & Q(E)=  \Xi (x,E)^2\left( E- \sum_{i=0}^3 l_i(l_i+1)\wp (x+\omega_i)\right) +\frac{1}{2}\Xi (x,E)\frac{d^2\Xi (x,E)}{dx^2}-\frac{1}{4}\left(\frac{d\Xi (x,E)}{dx} \right)^2. \label{const}
\end{align}
It is shown in \cite{Tak1} that $Q(E)$ is independent of $x$. Thus $Q(E)$ is a monic polynomial in $E$ of degree $2g+1$, which follows from the expression for $\Xi (x,E)$ given by (\ref{Fx}). The following proposition is obtained in \cite{Tak1}:
\begin{prop} \cite[Proposition 3.7]{Tak1}
Let $\Xi (x,E)$ be the doubly periodic function defined in Proposition \ref{prop:prod} and $Q(E)$ be the monic polynomial defined in (\ref{const}).
Then the function 
\begin{equation}
\Lambda (x,E)=\sqrt{\Xi (x,E)}\exp \int \frac{\sqrt{-Q(E)}dx}{\Xi (x,E)}
\label{eqn:Lam}
\end{equation}
is a solution to the differential equation (\ref{InoEF}).
\end{prop}

\subsection{The invariant subspace of elliptic functions} \label{sec:invsp}$ $

It is shown in \cite{Tak1} that there are finite-dimensional spaces of elliptic functions on which the action of the Hamiltonian (\ref{Ino}) is well defined. In this subsection, we introduce such spaces and discuss the relationship with the polynomial $Q(E)$ of the previous subsection.

Let ${\mathcal F}$ be the space spanned by meromorphic doubly periodic functions up to signs, namely
\begin{align}
& {\mathcal F}=\bigoplus _{\epsilon _1 , \epsilon _3 =\pm 1 } {\mathcal F} _{\epsilon _1 , \epsilon _3 }, \label{spaceF} \\
& {\mathcal F} _{\epsilon _1 , \epsilon _3 }=\{ f(x) \mbox{: meromorphic }| f(x+2\omega_1)= \epsilon _1 f(x), \; f(x+2\omega_3)= \epsilon _3 f(x) \} ,
\end{align}
where $(2\omega_1, 2\omega_3)$ are basic periods of elliptic functions.
%
 
Let $l_i \in \Zint _{\geq 0}$ $(i=0,1,2,3)$ and let $k_i$ be the rearrangement of $l_i$ such that $k_0\geq k_1 \geq k_2 \geq k_3 (\geq 0)$.
The dimension of the finite-dimensional invariant subspaces of the space ${\mathcal F}$ and the construction of the spaces are given in \cite[Theorem 3.1]{Tak1}. 

\begin{prop} \cite[Theorem 3.1]{Tak1} \label{periodic}

(i) If the number $k_0+ k_1+ k_2+ k_3 (=l_0+l_1+l_2+l_3)$ is even and $k_0+k_3\geq k_1+k_2$, then the dimension of the maximal finite-dimensional invariant subspace in ${\mathcal F}$ with respect to the action of the Hamiltonian $H$ (see (\ref{Ino})) is $2k_0+1$.

(ii) If the number $k_0+ k_1+ k_2+ k_3$ is even and $k_0+k_3< k_1+k_2$, then the dimension is $k_0+k_1+k_2-k_3+1$.

(iii) If the number $k_0+ k_1+ k_2+ k_3$ is odd and $k_0\geq k_1+k_2+k_3+1 $, then the dimension is $2k_0+1$.

(iv) If the number $k_0+ k_1+ k_2+ k_3$ is odd and $k_0< k_1+k_2+k_3+1$, then the dimension is $k_0+k_1+k_2+k_3+2$.
\end{prop}

Let $V$ be the finite-dimensional invariant space described in Proposition \ref{periodic} and $P(E)$ be the monic characteristic polynomial of the Hamiltonian $H$ (\ref{Ino}) on the space $V$, i.e. $P(E)=\prod_{i=1}^{\dim V} (E-E_i)$ ($\{ E_i\} $ are eigenvalues of $H$ on $V$).

Here we present the invariant space $V$ more explicitly. We assume $\tilde{\beta} _i\in \Zint$ $(i=0,1,2,3)$ and $-\sum_{i=0}^3 \tilde{\beta}_i /2 \in \Zint_{\geq 0}$. Let $V_{\tilde{\beta}_0, \tilde{\beta}_1, \tilde{\beta}_2, \tilde{\beta}_3}$ be the $(d+1)$-dimensional vector space spanned by  $\left\{ \wp_1(x)^{\tilde{\beta }_1}\wp_2(x)^{\tilde{\beta }_2}\wp_3(x)^{\tilde{\beta }_3}\wp(x)^n\right\} _{n=0, \dots ,d}$, where $d=-\sum_{i=0}^3 \tilde{\beta}_i /2$ and $\wp_i(x)$ ($i=1,2,3$) are co-$\wp$ functions (see appendix). 
We assume $\tilde{\alpha}_i \in \{ -l_i, l_i+1 \}$ $(i=0,1,2,3)$ and $\sum_{i=0}^3 \tilde{\alpha}_i/2 \in \Zint$. 
Then we set
\begin{align}
& U_{\tilde{\alpha}_0, \tilde{\alpha}_1, \tilde{\alpha}_2, \tilde{\alpha}_3}=
\left\{
\begin{array}{ll}
V_{\tilde{\alpha}_0, \tilde{\alpha}_1, \tilde{\alpha}_2, \tilde{\alpha}_3}, & \sum_{i=0}^3 \tilde{\alpha}_i/2 \in \Zint_{\leq 0}; \\
\{ 0 \} , & \sum_{i=0}^3 \tilde{\alpha}_i /2=1 ;\\
V_{1-\tilde{\alpha}_0, 1-\tilde{\alpha}_1, 1-\tilde{\alpha}_2, 1-\tilde{\alpha}_3}, & \sum_{i=0}^3 \tilde{\alpha}_i /2\in \Zint_{\geq 2}. 
\end{array}
\right. 
\end{align}

If $l_0 +l_1 +l_2 +l_3$ is even, then the space $V$ is decomposed as $V=U_{-l_0,-l_1,-l_2,-l_3}\oplus U_{-l_0 ,-l_1,l_2+1 ,l_3+1}\oplus U_{-l_0,l_1+1,-l_2,l_3+1}\oplus U_{-l_0 ,l_1+1,l_2+1,-l_3}$. 
If  $l_0 +l_1 +l_2 +l_3$ is odd, then we have $V=U_{-l_0,-l_1,-l_2,l_3+1}\oplus U_{-l_0 ,-l_1,l_2+1,-l_3} \oplus U_{-l_0 ,l_1+1,-l_2,-l_3} \oplus U_{l_0+1,-l_1,-l_2,-l_3}$.

For the equations $P(E)=0$ and $Q(E)=0$ we obtained the following proposition in \cite{Tak1}.

\begin{prop} (c.f. \cite[Theorem 3.8]{Tak1}) \label{prop:zeros}
The set of zeros of $Q(E)$ coincides with the set of zeros of $P(E)$.
\end{prop}

We propose the following conjecture.
\begin{conj} \label{conj1}
$P(E)=Q(E)$.
\end{conj}

By applying Proposition \ref{prop:zeros}, we immediately have the following proposition.

\begin{prop} \label{prop:Conject}

(i) If there are no multiple roots for each of the equations $P(E)=0$ and $Q(E)=0$ for generic $\tau $, then Conjecture \ref{conj1} is true.

(ii) If there are no multiple roots for the equation $Q(E)=0$ for generic $\tau$ and $\deg_E P(E)=\deg_E Q(E)$, then Conjecture \ref{conj1} is true.
\end{prop}

The following proposition follows upon combining Theorem 3.2 and Proposition 3.9 in \cite{Tak1} with Proposition \ref{prop:Conject} (ii) of this paper. 
\begin{prop} \label{propconj1}
If two of the $l_i$ ($i=0,1,2,3$) are zero, then Conjecture \ref{conj1} is true.
\end{prop}

If $l_0+l_1+l_2+l_3$ is small, then we can confirm  Conjecture \ref{conj1} directly through direct calculation.
\begin{prop}  (c.f. \cite[Proposition 3.10]{Tak1})
If  $l_0+l_1+l_2+l_3\leq 8$, then Conjecture \ref{conj1} is true.
\end{prop}

\section{Finite-gap property} \label{sec3}

\subsection{The construction of the commuting operator}$ $ \label{sec:constr}

Let us consider the Schr\"odinger operator $-d^2/dx^2+q(x)$ with the conditions that $q(x)$ is meromorphic on $\Rea$ (may have poles) and $q(x+1)=q(x)$.
These conditions on $q(x)$ are assumed in the definitions below.

\begin{df} \label{def:fin}
Assume that $q(x) \in \Rea $ for $x \in \Rea$ and there are no poles for $q(x)$ on $\Rea$. We set $H=-d^2/dx^2+q(x)$. Let $\sigma _b(H)$ be the set such that 
$$
E \in \sigma _b(H) \; \Leftrightarrow \mbox{ All solutions to }(H-E)f(x)=0 \mbox{ are bounded on }x \in \Rea 
$$
and $\overline{\sigma _b(H)}$ be the closure of the set $\sigma _b(H)$. If the set $\Rea \setminus \overline{\sigma _b(H)}$ can be written as
\begin{equation}
\Rea \setminus \overline{\sigma _b(H)}= (-\infty,E_{0})\cup (E_{1},E_{2})\cup \dots \cup (E_{2g-1}, E_{2g})
\end{equation}
with $E_0<E_{1}<\cdots <E_{2g}$, then $q(x)$ is called the finite-gap potential.
\end{df}
%

\begin{df} \label{def:algfin}
If there exists an odd-order differential operator 
$A= \left( d/dx \right)^{2g+1} +  \! $ $ \sum_{j=0}^{2g-1}\! $ $ b_j(x) \left( d/dx \right)^{2g-1-j} $ such that $[A, -d^2/dx^2+q(x)]=0$, then $q(x)$ is called the algebro-geometric finite-gap potential.
\end{df}

\begin{remk}
It is known that the equation  $[A, -d^2/dx^2+q(x)]=0$ is equivalent to the function $q(x)$ being a solution to some stationary higher-order KdV equation.
\end{remk}

It is known that, if $q(x) $ is real-holomorphic on $\Rea$ and $q(x+1)=q(x)$, then $q(x)$ is a finite-gap potential if and only if $q(x)$ is an algebro-geometric finite-gap potential (for details see \cite{Nov,Dub,Wei}).

In this section it is shown that $\sum_{i=0}^3 l_i(l_i+1)\wp (x+\omega_i)$ is an algebro-geometric finite-gap potential if $l_i \in \Zint $ $(i=0,1,2,3)$.
By the substitution $l_i \leftrightarrow -l_i-1$, it is sufficient to show it for the case $l_i \in \Zint _{\geq 0}$ $(i=0,1,2,3)$.

Write
\begin{equation}
\Xi(x,E) = \sum_{i=0}^{g} a_{g-i}(x) E^i. \label{Xiag}
\end{equation}
From Proposition \ref{prop:prod} we have $a_0(x)=1$.

\begin{thm} \label{thm:PhiA}
Set $u(x)=  \sum_{i=0}^3 l_i(l_i+1)\wp (x+\omega_i)$ and define the $(2g+1)$st-order differential operator $A$ by
\begin{equation}
A= \sum_{j=0}^{g} \left\{ a_j(x)\frac{d}{dx}-\frac{1}{2} \left( \frac{d}{dx} a_j(x) \right) \right\} \left( - \frac{d^2}{dx^2} +u(x) \right) ^{g-j}, \label{Adef}
\end{equation}
where the $a_j(x)$ are defined in (\ref{Xiag}).
Then the operator $A$ commutes with the Hamiltonian $H=-d^2/dx^2  +u(x) $. In other words the function $u(x)=  \sum_{i=0}^3 l_i(l_i+1)\wp (x+\omega_i)$ is an algebro-geometric finite-gap potential.
\end{thm}
\begin{proof}
Since the function $\Xi (x,E)$ in (\ref{Xiag}) satisfies the differential equation (\ref{prodDE}), we obtain the following relations by equating the coefficients of $E^{g-j}$:
\begin{equation}
a'''_j(x)-4u(x)a'_j(x)-2u'(x)a_j(x)+4a'_{j+1}(x)=0.
\end{equation}

On the other hand
\begin{align}
& [A,H]= \sum_{j=0}^{g} \left[ a_j(x)\frac{d}{dx}-\frac{1}{2} a'_j(x), H\right] H^{g-j} \\
& =  \sum_{j=0}^{g} \left(-2a'_j(x)\frac{d^2}{dx^2}+\frac{1}{2} a'''_j(x) -u'(x)a_j(x)\right)H^{g-j} \nonumber \\
& = \sum_{j=0}^{g} \left(-2a'_j(x)(-H+u(x))+\frac{1}{2} a'''_j(x) -u'(x)a_j(x)\right)H^{g-j} \nonumber \\
& = \frac{1}{2} \sum_{j=0}^{g} \left(4a'_{j+1}(x)-4u(x)a'_j(x)+ a'''_j(x) -2u'(x)a_j(x)\right)H^{g-j} \nonumber \\
& =0. \nonumber
\end{align}
\end{proof}

\begin{prop} \label{prop:algrel}
Let $Q(E)$ be the polynomial defined in (\ref{const}). Then
\begin{equation} 
A^2+Q(H)=0. \label{algrel}
\end{equation}
\end{prop}
\begin{proof}
Let $f(x)$ be a solution to the second-order differential equation $(H-E)f(x)=0$.
From the definitions of $A$ in (\ref{Adef}) and (\ref{Xiag}) it follows that
\begin{align}
Af(x)& = \sum_{j=0}^{g} \left\{ a_j(x)E^{g-j} \frac{d}{dx}-\frac{1}{2} \left( \frac{d}{dx} a_j(x)E^{g-j} \right) \right\} f(x) \\
& = \Xi (x,E) f'(x) -\frac{1}{2}\Xi '(x,E) f(x). \nonumber
\end{align}
Here $\Xi '(x,E) := \frac{\partial}{\partial x}\Xi (x,E) $. From the commutativity of $A$ and $H$ we have $(H-E)Af(x)=0$. Hence
\begin{align}
A(Af(x))& =  \Xi (x,E)\frac{d}{dx}\left(A f(x)\right)  -\frac{1}{2}\Xi '(x,E) A f(x) \\
& = \Xi (x,E)^2 f''(x) -\frac{1}{2}\Xi ''(x,E) \Xi (x,E)f(x)+\frac{1}{4}\Xi '(x,E)^2 f(x). \nonumber \\
& = \left( \Xi (x,E)^2 (u(x)-E)-\frac{1}{2}\Xi ''(x,E)\Xi (x,E)+\frac{1}{4}\Xi '(x,E)^2 \right) f(x) \nonumber \\
& =- Q(E) f(x)= -Q(H)f(x),\nonumber 
\end{align}
where relation (\ref{const}) was used.
The operator $A^2 +Q(H)$ is a differential operator with degree at most $2g+1$.
Since the equality $(A^2+Q(H))f(x)=0$ is valid for all eigenfunctions of $H$, it follows that $A^2+Q(H)=0$.
\end{proof}

\begin{remk}
Theorem \ref{thm:PhiA} and Proposition \ref{prop:algrel} are established in \cite{IM,TD} for a smooth finite-gap potential.
\end{remk}

Let $\Gamma \! : y^2=-Q(x)$ be the spectral curve of the operators $A$ and $H$. This curve plays an important role in the formulas (\ref{eqn:Lam}) and (\ref{algrel}). By Proposition \ref{prop:zeros} edges of the hyperelliptic curve $\Gamma $ are eigenvalues on the invariant space $V$ (see section \ref{sec:invsp}). 
The genus of the curve $\Gamma $ is $g$, where $g$ is as in (\ref{Xiag}).

\begin{prop} \label{prop:genus}
Let $k_0= \max (l_0, l_1, l_2, l_3)$, $k_3=\min (l_0, l_1, l_2, l_3)$ and $l=l_0+l_1+l_2+l_3$. If Conjecture \ref{conj1} is true, then the arithmetic genus of the curve $\Gamma $ is given by 
\begin{equation}
\begin{array}{ll}
k_0, & l \mbox{ is even and } k_0+k_3\geq l/2 ;\\
( l-2k_3 )/2 , & l \mbox{ is even and } k_0+k_3< l/2 ;\\
k_0, & l \mbox{ is odd and } k_0\geq ( l +1)/2 ;\\
( l+1 )/2 , & l \mbox{ is odd and } k_0< ( l +1 )/2 .
\end{array}
\label{genus}
\end{equation}
\end{prop}
\begin{proof}
It follows directly from Proposition \ref{periodic}.
\end{proof}

\begin{remk}
The values in (\ref{genus}) have already appeared in \cite{GW,Smi,Tak1}, although the present author could not understand the proofs of computing the genus of $\Gamma $ given in \cite{GW,Smi} because he could not confirm how to avoid the possibility that either of the equations $P(E)=0$ and $Q(E)=0$ could have multiple roots.
\end{remk}

By combining Propositions \ref{propconj1} and \ref{prop:genus} we obtain the following corollary. 
\begin{cor}
If two of $l_i$ ($i=0,1,2,3$) are zero, then the arithmetic genus of the curve $\Gamma $ is $\max (l_0, l_1, l_2, l_3)$.
\end{cor}
%
\begin{prop} \label{prop:fingapspect}
If $l_0= l_1 =0$, $l_2 ,l_3 \in \Zint_{\geq 0}$ and $\tau \in \Zint + \sqrt{-1} \Rea _{>0}$, then $u(x)=  \sum_{i=2}^3 l_i(l_i+1)\wp (x+\omega_i)$ is a finite-gap potential.
\end{prop}
\begin{proof}
First from the condition $\tau \in \Zint + \sqrt{-1} \Rea _{>0}$ we have $u(x) \in \Rea $ for $x\in \Rea$.
Also, since the poles of the function $\wp (x+ \tau /2)$ (resp. $\wp (x+ (\tau +1)/2)$) are located on the set $\Zint +(\Zint +1/2)\tau$ (resp. $(\Zint+1/2) +(\Zint +1/2)\tau$ ), the function $u(x)$ is holomorphic on $\Rea$.
Finally it is shown in \cite{Nov,Dub} that, if $u(x) $ is real-holomorphic on $\Rea$ and $u(x+1)=u(x)$, then $u(x)$ is a finite-gap potential if and only if $u(x)$ is an algebro-geometric finite-gap potential.
Thus we obtain the proposition.
\end{proof}

\subsection{The determinantal formula for the commuting operator}$ $ \label{sec:det}

In this subsection a determinantal formula for the commuting operator $A$ is introduced.

\begin{conj} \label{conj:conj2}
Let $l_i$ ($i=0,1,2,3$) be nonnegative integers and $V$ be the invariant subspace defined in section \ref{sec:invsp}.

Let $f_1(x), \dots , f_{2\tilde{g}+1}(x)$ be a basis for the space $V$. Then
\begin{equation}
A=A_0 \left| \left(
\begin{array}{cccc}
f_1(x) & \dots & f_{2\tilde{g}+1}(x) & 1 \\
\frac{d}{dx}f_1(x) & \dots &\frac{d}{dx} f_{2\tilde{g}+1}(x) & \frac{d}{dx} \\
\left(\frac{d}{dx} \right)^2 f_1(x) & \dots &\left(\frac{d}{dx} \right)^2 f_{2\tilde{g}+1}(x) & \left(\frac{d}{dx} \right)^2 \\
\vdots & & \vdots & \vdots \\
\left(\frac{d}{dx} \right)^{2\tilde{g}+1} f_1(x) &  \dots &\left(\frac{d}{dx} \right)^{2\tilde{g}+1} f_{2\tilde{g}+1}(x) & \left(\frac{d}{dx} \right)^{2\tilde{g}+1}
\end{array} \right) \right|  \label{conj2rel}
\end{equation}
for some nonzero constant $A_0$. Here the determinant for the matrix 
$ \left( \begin{array}{ccc}
a_{1,1} &  \dots &  a_{1,n} \\
\vdots  &  & \vdots \\
a_{n,1} &  \dots &  a_{n,n}
\end{array} \right) $
is defined by $\sum _{\sigma \in {\frak S}_n}\mbox{sgn}(\sigma) a_{\sigma(1),1} \dots  a_{\sigma(n),n}$.
\end{conj}

\begin{prop} \label{prop:propconj2}
If two of $l_i$ ($i=0,1,2,3$) are zero, Conjecture \ref{conj:conj2} is true.
\end{prop}
\begin{proof}
It is sufficient to show the statements for the case $l_2=l_3=0$ and $l_0\geq l_1$. From Propositions \ref{propconj1} and \ref{prop:algrel} we have $A^2=-P(H)$.
In \cite[Theorem 3.2]{Tak1} it is shown that the equation $P(E)=0$ has no multiple roots for generic $\tau$.
For this case we have $\dim V=2g+1$, where $g$ is defined in Proposition \ref{prop:prod}.

We fix the period $\tau _0$ to satisfy the condition  that the equation $P(E)=0$ has no multiple roots.
Let $E_i$ $(i=1,\dots ,2g+1)$ be the zeros of $P(E)=0$ and let $\tilde{V}_i$ be the $2$-dimensional space spanned by solutions to the differential equation $(H-E_i)f(x)=0$.

From the commutativity $[A,H]=0$ the operator $A$ preserves the space $\tilde{V} _i$. For $f(x) \in \tilde{V} _i$ we have $A^2 f(x)= -P(H)f(x)=-P(E_i)f(x)=0$. Therefore all eigenvalues of $A$ on $\tilde{V} _i$ are $0$.

Let $\xi _i(x) \in V \cap \tilde{V} _i$. From the definition of $V$ the function $\xi _i(x) $ is doubly periodic up to signs and satisfies the equation $(H-E_i)\xi _i(x)=0$. 
Another solution to $(H-E_i) f(x)=0$ is given by 
$\xi _i(x) \int _0^x  dy/\xi _i(y) ^2$. This function could never be a doubly periodic function up to signs because it has $l_0 +1$ zeros at $x=0$ and has fewer poles on the fundamental domain. Hence the dimension of the space of doubly periodic functions up to signs satisfying the differential equation $(H-E_i)f(x)=0$ is one.
Since the function $\Xi(x,E)$ is doubly periodic, the coefficients of $\left( d/dx\right)^j$ $(j=0, \dots ,2g+1)$ on $A$ are doubly periodic (see (\ref{Adef})); hence the function $A \xi _i(x)$ is proportional to $\xi _i(x)$. Thus we have $A \xi _i(x)=0$, because the eigenvalue must be zero and
therefore we also have $Af(x)=0$ for all $f(x) \in V$.

Since $\deg A=\dim V$, all solutions to $A\tilde{f}(x)=0$ must satisfy $\tilde{f}(x) \in V$ and we have 
\begin{equation}
A= h(x) \left| \left(
\begin{array}{cccc}
f_1(x) & \dots & f_{2g+1}(x) & 1 \\
\frac{d}{dx}f_1(x) & \dots &\frac{d}{dx} f_{2g+1}(x) & \frac{d}{dx} \\
\left(\frac{d}{dx} \right)^2 f_1(x) & \dots &\left(\frac{d}{dx} \right)^2 f_{2g+1}(x) & \left(\frac{d}{dx} \right)^2 \\
\vdots & & \vdots & \vdots \\
\left(\frac{d}{dx} \right)^{2g+1} f_1(x) &  \dots &\left(\frac{d}{dx} \right)^{2g+1} f_{2g+1}(x) & \left(\frac{d}{dx} \right)^{2g+1}
\end{array} \right) \right|  \label{Agx}
\end{equation}
for some function $h(x)$. 
Since the coefficient of $ \left(d/dx \right)^{2g+1}$ on the right-hand side of (\ref{Agx}) is $(-1)^g$ (see (\ref{Adef})), we have
\begin{equation}
\left| \left(
\begin{array}{ccc}
f_1(x) & \dots & f_{2g+1}(x)  \\
\frac{d}{dx}f_1(x) & \dots &\frac{d}{dx} f_{2g+1}(x) \\
\vdots & & \vdots \\
\left(\frac{d}{dx} \right)^{2g} f_1(x) &  \dots &\left(\frac{d}{dx} \right)^{2g} f_{2g+1}(x) 
\end{array} \right) \right|
= \frac{(-1)^g}{h(x)}. \label{gx-1}
\end{equation}
By differentiating equality (\ref{gx-1}), we obtain that the coefficient of $ \left(d/dx \right)^{2g}$ on the right-hand side of (\ref{Agx}) is equal to $(-1)^g h'(x)/h(x)$.
From the commutativity $[A,H]=0$ the coefficient of $\left( d/dx \right) ^{2g} $ in the operator $A$ must be zero, whence $h'(x)=0$.

Therefore relation (\ref{conj2rel}) follows for the case $l_2=l_3=0$ and $l_0\geq l_1$, if the equation $P(E)=0$ has no multiple roots. From analytic continuation in $\tau$ we obtain relation (\ref{conj2rel}) for all $\tau \in \Rea +\sqrt{-1}\Rea _{>0}$ when $l_2=l_3=0$, $l_0\geq l_1$.
\end{proof}

\subsection{Monodromy}$ $ \label{sec:monod}

In this subsection we calculate the monodromy, specifically the constant $B(E)$ for $\Lambda (x+1,E)= B(E)\Lambda(x,E)$, where the function $\Lambda (x,E)$ is the solution to the equation $(H-E)\Lambda(x,E) =0$ of the form (\ref{eqn:Lam}).
The monodromy plays an important role in investigating eigenvalues and eigenstates for systems with physical boundary conditions.

From formula (\ref{eqn:Lam}) we have 
\begin{equation}
\Lambda (x+1,E)=\pm \Lambda (x,E) \exp \int _{0+ \varepsilon }^{1+\varepsilon }\frac{\sqrt{-Q(E)}dx}{\Xi (x,E)},
\label{eqn:Lam01}
\end{equation}
with $\varepsilon $ a constant so determined as to avoid passing through the poles in the integration. The sign $\pm$ is determined by the analytic continuation of the function $\sqrt{\Xi (x,E)}$.

Although this formula is explicit, the calculation of the integral would be difficult.

In \cite{Tak1} it is shown that the function $\Lambda (x,E)$ admits the following expression of the Bethe ansatz form:
\begin{equation}
\Lambda (x,E)=\frac{\exp(\pi \sqrt{-1}cx)\prod_{j=1}^l \theta(x+t_j)}{\theta(x)^{l_0}\theta(x+1/2)^{l_1}\theta(x+(1+\tau)/2)^{l_2}\theta(x+\tau/2)^{l_3}},
\label{tBV}
\end{equation}
where $l=l_0+l_1+l_2+l_3$ and $\theta (x)$ is the elliptic theta function.
The values $t_j$ $(j=1,\dots ,l)$ and $c$ satisfy the Bethe ansatz equation given in \cite[Theorem 3.13]{Tak1}.

From this expression we have
\begin{equation}
\Lambda (x+1,E)=\Lambda (x,E) \exp (\pi \sqrt{-1}c).
\end{equation}
To obtain the value $c$ we must solve the Bethe ansatz equation, which is in general transcendental.

In this subsection we introduce a different formula. The monodromy is expressed as an hyperelliptic integral of second kind.
The following formula for a smooth finite-gap potential is found in \cite{Hoc,MvM}.

\begin{thm} \label{thm:conj3}
Let $l_i \in \Zint_{\geq 0}$, $P \in 2 \Zint \omega_1 \oplus 2 \Zint \omega_3 (=\Zint \oplus \Zint \tau)$, and $E^{(P)}_0$ be the eigenvalue of $H$ such that $\Lambda (x+P,E^{(P)}_0)=\Lambda (x,E^{(P)}_0)$ $($resp. $\Lambda (x+P,E^{(P)}_0)=-\Lambda (x,E^{(P)}_0))$. Then
\begin{align} 
& \Lambda (x+P,E_*)=\Lambda (x,E_*) \exp \left( -\frac{1}{2} \int_{E^{(P)}_0}^{E_*}\frac{ \int_{0+\varepsilon }^{P+\varepsilon }\Xi (x,E)dx}{\sqrt{-Q(E)}} dE\right)
\label{analcontP}\\
& \left( \mbox{resp. } \Lambda (x+P,E_*)=-\Lambda (x,E_*) \exp \left( -\frac{1}{2} \int_{E^{(P)}_0}^{E_*}\frac{ \int_{0+\varepsilon }^{P+\varepsilon }\Xi (x,E)dx}{\sqrt{-Q(E)}} dE\right) \right) \nonumber 
\end{align}
with $\varepsilon $ a constant so determined as to avoid passing through the poles in the integration. 
In particular, if $E_0$ is the eigenvalue of $H$ such that $\Lambda (x+1,E_0)=\Lambda (x,E_0)$ $($resp. $\Lambda (x+1,E_0)=-\Lambda (x,E_0))$, then
\begin{align} 
& \Lambda (x+1,E_*)=\Lambda (x,E_*) \exp \left( -\frac{1}{2} \int_{E_0}^{E_*}\frac{Q_1(E)}{\sqrt{-Q(E)}} dE\right)
\label{analcont} \\
& \left( \mbox{resp. } \Lambda (x+1,E_*)=-\Lambda (x,E_*) \exp \left( -\frac{1}{2} \int_{E_0}^{E_*}\frac{Q_1(E)}{\sqrt{-Q(E)}} dE\right) \right), \nonumber
\end{align}
where the polynomial $Q_1(E)$ of degree $g$ is given by 
\begin{equation}
Q_1(E)= \int_{0+\varepsilon }^{1+\varepsilon }\Xi (x,E)dx . \label{analcontnum}
\end{equation}
\end{thm}
\begin{proof}
Set 
\begin{equation}
R=R(x,E)=\frac{\Xi(x,E)}{\sqrt{-Q(E)}}.
\end{equation}
We use $'$ instead of $\partial/\partial x$. Set $v=\sum_{i=0}^3 l_i(l_i+1)\wp (x+\omega_i)$.
From (\ref{const})
\begin{align}
 & R ^2\left( v- E \right) -\frac{1}{2}RR''+\frac{1}{4}\left(R' \right)^2=1. \label{const1}
\end{align}
By differentiating (\ref{const1}) with respect to $E$ we have
\begin{equation}
-\frac{\partial R''}{\partial E}=-\frac{R'}{R}\frac{\partial R'}{\partial E}+2R+\frac{R''}{R}\frac{\partial R}{\partial E}-4(v-E)\frac{\partial R}{\partial E}.
\end{equation}
Then
\begin{align}
& \left( -\frac{\partial R'}{\partial E} +\frac{R'}{R}\frac{\partial R}{\partial E}\right) ' \label{eq:R} \\
& =\left(-\frac{R'}{R}\frac{\partial R'}{\partial E}+2R+\frac{R''}{R}\frac{\partial R}{\partial E}-4(v-E)\frac{\partial R}{\partial E}\right)+ \left( \frac{R''}{R}-\frac{(R')^2}{R^2} \right) \frac{\partial R}{\partial E}+\frac{R'}{R}\frac{\partial R'}{\partial E}\nonumber \\
& =2R+\frac{1}{R^2}\frac{\partial R}{\partial E}\left( 2RR''-4R^2(v-E)-(R')^2\right) \nonumber \\
& = 2R-\frac{4}{R^2}\frac{\partial R}{\partial E}=2R+4\frac{\partial }{\partial E}\left( \frac{1}{R} \right). \nonumber 
\end{align}
Since the function $-\frac{\partial R'}{\partial E} +\frac{R'}{R}\frac{\partial R}{\partial E}$ is doubly periodic in $x$, it follows that 
\begin{equation}
\int_{0+\varepsilon }^{P+\varepsilon } \left( -\frac{\partial R'}{\partial E} +\frac{R'}{R}\frac{\partial R}{\partial E}\right) 'dx =0.
\end{equation}
Hence
\begin{equation}
-\frac{1}{2} \int_{0+\varepsilon }^{P+\varepsilon } Rdx=\int_{0+\varepsilon }^{P+\varepsilon }\frac{\partial }{\partial E}\left( \frac{1}{R} \right)dx.
\end{equation}
By integrating with respect to $E$ we see that
\begin{equation}
\exp \left( -\frac{1}{2} \int_{E_0^{(P)}}^{E_*}\frac{\int_{0+\varepsilon }^{P+\varepsilon } \Xi(x,E)dx}{\sqrt{-Q(E)}}dE \right) = \frac{\exp \left( \int_{0+\varepsilon }^{P+\varepsilon }\frac{\sqrt{-Q(E)}}{\Xi(x,E)} dx \right)}{\exp \left( \int_{0+\varepsilon }^{P+\varepsilon }\frac{\sqrt{-Q(E^{(P)}_0)}}{\Xi(x,E^{(P)}_0)} dx \right)} ,
\label{inteq}
\end{equation}

From formula (\ref{eqn:Lam}) we have 
\begin{equation}
\Lambda (x+P,E)=\pm \Lambda (x,E) \exp \int _{0+ \varepsilon }^{P+\varepsilon }\frac{\sqrt{-Q(E)}dx}{\Xi (x,E)}.
\label{eqn:Lam01P}
\end{equation}
The sign $\pm$ is determined by the analytic continuation of the function $\sqrt{\Xi (x,E)}$.
By the assumption $\Lambda (x+P,E^{(P)}_0)=\Lambda (x,E^{(P)}_0)$ (resp. $\Lambda (x+P,E^{(P)}_0)=-\Lambda (x,E^{(P)}_0)$) we have $\pm \exp \int _{0+ \varepsilon }^{P+\varepsilon }\frac{\sqrt{-Q(E^{(P)}_0)}dx}{\Xi (x,E^{(P)}_0)} = 1$ (resp. $\pm \exp \int _{0+ \varepsilon }^{P+\varepsilon }\frac{\sqrt{-Q(E^{(P)}_0)}dx}{\Xi (x,E^{(P)}_0)} = -1$), where the sign $\pm$ is chosen to be the same as the one in (\ref{eqn:Lam01P}).
Combining with (\ref{inteq}) and (\ref{eqn:Lam01P}) we obtain (\ref{analcontP}).
\end{proof}

\begin{remk}
(i) A similar equality to (\ref{eq:R}) can be found in \cite{GD}.\\
(ii) Although the integral $\int_{E_0}^{E_*}\frac{Q_1(E)}{\sqrt{-Q(E)}} dE$ may be dependent on the choice of the initial point $E_0$ such that $\Lambda (x+1,E_0) = \Lambda (x,E_0)$ (resp. $\Lambda (x+1,E_0) = -\Lambda (x,E_0)$) and the path from $E_0$ to $E_*$, the function $\exp \left( -\frac{1}{2}\int_{E_0}^{E_*}\frac{Q_1(E)}{\sqrt{-Q(E)}} dE\right)+ \exp \left( \frac{1}{2}\int_{E_0}^{E_*}\frac{Q_1(E)}{\sqrt{-Q(E)}} dE\right)$ is independent on $E_0$ with $\Lambda (x+1,E_0) = \Lambda (x,E_0)$ (resp. $\Lambda (x+1,E_0) = -\Lambda (x,E_0)$) and the path.
\end{remk}

From Proposition \ref{prop:zeros}, if $Q(E_0)=0$, then the differential equation $(H-E_0)f(x)=0$ (see (\ref{InoEF})) has a solution of doubly periodic function up to signs. It is seen that the solution is written as $\Lambda (x,E_0)=\sqrt{\Xi(x,E_0)}$ up to constant multiplication. Hence $\Lambda (x+1,E_0) = \Lambda (x,E_0)$ or $\Lambda (x+1,E_0) = -\Lambda (x,E_0)$.
Combined with Theorem \ref{thm:conj3}, we obtain the following corollary:

\begin{cor} \label{cor:conj3}
Let $l_i \in \Zint_{\geq 0}$ and let $E_0$ be a solution to $Q(E)=0$. Then  $\Lambda (x+1,E_0) = \Lambda (x,E_0)$ or $\Lambda (x+1,E_0) = -\Lambda (x,E_0)$.
If $\Lambda (x+1,E_0) = \Lambda (x,E_0)$, then $\Lambda (x+1,E_*)=\Lambda (x,E_*) \exp \left( -\frac{1}{2} \int_{E_0}^{E_*}\frac{Q_1(E)}{\sqrt{-Q(E)}} dE\right)$ and, if $\Lambda (x+1,E_0) = -\Lambda (x,E_0)$, then $\Lambda (x+1,E_*)=-\Lambda (x,E_*) \exp \left( -\frac{1}{2} \int_{E_0}^{E_*}\frac{Q_1(E)}{\sqrt{-Q(E)}} dE\right)$.
\end{cor}

Theoretically the integral $Q_1(E)=\int_{0+\varepsilon }^{1+\varepsilon }\Xi (x,E)dx$ can be calculated explicitly. This is illustrated by some examples in the next subsection.

\subsection{Examples} $ $ \label{sec:exam}

Two examples are provided to illustrate plainly what has been done to this point.

\subsubsection{The case $l_0=0, l_1=0, l_2=0, l_3=1$} $ $

The Hamiltonian $H$ is given by
\begin{equation}
H= -\frac{d^2}{dx^2} + 2\wp \left(x+\frac{\tau}{2}\right).
\label{Ino0001}
\end{equation}
For this case the space of finite-dimensional elliptic functions $V$ defined in section \ref{sec:invsp} is spanned by co-$\wp $ functions $\wp_1(x+\frac{\tau}{2})$, $\wp_2(x+\frac{\tau}{2})$ and $\wp_3(x+\frac{\tau}{2})$. Hence $\dim V=3$.
The characteristic polynomial of the Hamiltonian $H$ on the space $V$ is
\begin{equation}
P(E)=(E+e_1)(E+e_2)(E+e_3),
\end{equation}
where $\wp(\frac{1}{2})=e_1$,  $\wp(\frac{1+\tau}{2})=e_2$, and $\wp(\frac{\tau}{2})=e_3$. For this case we have $Q(E)=P(E)$ from Proposition \ref{prop:propconj2}.

From Proposition \ref{prop:fingapspect} $ 2\wp (x+\frac{\tau}{2})$ is a finite-gap potential if $\tau \in \Zint + \sqrt{-1} \Rea _{>0}$. Indeed, if $\tau \in 2 \Zint + \sqrt{-1} \Rea _{>0}$, then we have
\begin{equation}
\Rea \setminus \overline{\sigma _b(H)}= (-\infty,-e_1) \cup (-e_2,-e_3),
\end{equation}
for  $-e_1<-e_2<-e_3$ (see Definition \ref{def:fin}).
If $\tau \in 2 \Zint + 1+ \sqrt{-1} \Rea _{>0}$, then $-e_1<-e_3<-e_2$ and 
\begin{equation}
\Rea \setminus \overline{\sigma _b(H)}= (-\infty,-e_1) \cup (-e_3,-e_2).
\end{equation}

The functions $\Xi(x,E)$ and $\Lambda (x,E)$ are given by
\begin{align}
& \Xi(x,E)=\wp \left(x+\frac{\tau}{2} \right)+E, \\
& \Lambda (x,E)= \left( \sqrt{\wp \left(x+\frac{\tau}{2}\right)+E} \right) \exp\int^x \frac{\sqrt{-(E+e_1)(E+e_2)(E+e_3)}}{\wp (x+\frac{\tau}{2})+E} dx. \nonumber
\end{align}
Here $\Lambda (\pm x,E)$ satisfy the equations $(H-E)\Lambda (\pm x,E)=0$ and $\Xi(x,E)$ satisfies equation (\ref{prodDE}).

The function $ 2\wp (x+\frac{\tau}{2})$ is also an algebro-geometric potential if $\tau \in \Rea +  \sqrt{-1} \Rea _{>0}$. The commuting operator $A$ (\ref{Adef}) is given by 
\begin{equation}
A=-\left( \left( \frac{d}{dx} \right) ^3-3 \wp\left(x+\frac{\tau}{2}\right)\frac{d}{dx}-\frac{3}{2} \wp '\left(x+\frac{\tau}{2}\right) \right) .
\end{equation}
As in the case of Proposition \ref{prop:propconj2}, we obtain the following determinantal formula:
\begin{align}
& (e_1-e_2)(e_2-e_3)(e_1-e_3)A \\
& \; \; \; =\left| \left(
\begin{array}{cccc}
\wp_1(x+\frac{\tau}{2})& \wp_2(x+\frac{\tau}{2})& \wp_3(x+\frac{\tau}{2})& 1 \\
\wp'_1(x+\frac{\tau}{2})& \wp'_2(x+\frac{\tau}{2})& \wp'_3(x+\frac{\tau}{2})& \frac{d}{dx} \\
\wp''_1(x+\frac{\tau}{2})& \wp''_2(x+\frac{\tau}{2})& \wp''_3(x+\frac{\tau}{2})& \left(\frac{d}{dx} \right)^2 \\
\wp'''_1(x+\frac{\tau}{2})& \wp'''_2(x+\frac{\tau}{2})& \wp'''_3(x+\frac{\tau}{2})& \left(\frac{d}{dx} \right)^3 
\end{array} \right) \right| . \nonumber
\end{align}
For the operators $A$ and $H$ there is an algebraic relation
\begin{equation}
A^2=-(H+e_1)(H+e_2)(H+e_3)
\end{equation}
and the genus of this curve is one.

Now we can proceed to calculate the coefficient of the monodromy (\ref{analcont}).
We have
\begin{align}
& \int_{0+\varepsilon }^{1+\varepsilon }\Xi (x,E)dx= 
\int_{0+\varepsilon }^{1+\varepsilon } \left( E+\wp \left(x+\frac{\tau}{2} \right) \right) dx \\
& \; \; \; \; =E-\zeta \left(1+\frac{\tau}{2}+\varepsilon \right)+\zeta\left(\frac{\tau}{2}+\varepsilon \right)=
E-2\eta_1 ,\nonumber
\end{align}
where $\eta_1 = \zeta (\frac{1}{2})$. (Here we used definitions and formulas from the appendix.)
Since $E=-e_1$ is a solution to $Q(E)=0$ such that $\Lambda (x+1,E_0)=\Lambda (x,E_0)$, it follows that 
\begin{equation} 
\Lambda (x+1,E_*)=\Lambda (x,E_*) \exp \left( -\frac{1}{2}\int_{-e_1}^{E_*}\frac{E-2\eta_1}{\sqrt{-(E+e_1)(E+e_2)(E+e_3)}} dE\right),
\label{analcont0001}
\end{equation}
which is an elliptic integral of second kind.

\subsubsection{The case $l_0=2, l_1=1, l_2=0, l_3=0$} $ $

The Hamiltonian $H$ is given by
\begin{equation}
H= -\frac{d^2}{dx^2} + 6\wp(x)+2\wp \left(x+\frac{1}{2}\right).
\label{Ino2100}
\end{equation}
For this case the space of finite-dimensional elliptic functions $V$ defined in section \ref{sec:invsp} is spanned by the functions $\wp(x) -e_1$, $\wp_2(x)/\wp_1(x)$, $\wp_3(x)/\wp_1(x)$, $\wp_1(x)\wp_2(x)$, and $\wp_1(x)\wp_3(x)$, so $\dim V=5$.
The characteristic polynomial of the Hamiltonian $H$ on the space $V$ is
\begin{align}
& P(E)=\frac{1}{16}(E+4e_1)\{  16E^4-224e_1E^3+(56g_2-144e_1^2)E^2 \label{2100P} \\
& \; \; \; \; \; -(200g_2e_1-5920e_1^3)E+(81g_2^2-3640g_2e_1^2+19216e_1^4)\}, \nonumber
\end{align}
where $g_2=-4(e_1e_2+e_2e_3+e_3e_1)$.  For this case we have $Q(E)=P(E)$ from Proposition \ref{prop:propconj2}.

The functions $\Xi(x,E)$ and $\Lambda (x,E)$ are given by
\begin{align}
 \Xi(x,E)=& 9\wp (x)^2 +3(E-2e_1)\wp(x)+(E+4e_1)\wp\left(x+\frac{1}{2}\right) +E^2-5e_1E-27e_1^2,  \\
 \Lambda (x,E)= & \sqrt{\Xi(x,E)}\exp\int^x \frac{\sqrt{-Q(E)}}{\Xi(x,E)} dx. \nonumber
\end{align}
Here $\Lambda (\pm x,E)$ satisfy the equations $(H-E)\Lambda (\pm x,E)=0$ and $\Xi(x,E)$ satisfies equation (\ref{prodDE}).

The function $ 6\wp(x)+2\wp (x+\frac{1}{2})$ is an algebro-geometric finite-gap potential. The commuting operator $A$ (\ref{Adef}) is given by 
\begin{align}
& A=\left( \frac{d}{dx} \right) ^5-5\left( 3\wp (x)+ \wp \left(x+\frac{1}{2}\right)-e_1 \right) \left(\frac{d}{dx}\right)^3 \\
& -\frac{15}{2} \left( 3\wp '(x) +\wp '\left(x+\frac{1}{2}\right) \right) \left(\frac{d}{dx}\right)^2  \nonumber \\
& -15\left( 3\wp(x)^2 +2 \wp \left(x+\frac{1}{2}\right) ^2 -2e_1 \wp \left(x+\frac{1}{2}\right) -3e_1^2 -\frac{1}{5}g_2 \right) \frac{d}{dx} \nonumber \\
& -15\left( \wp \left(x+\frac{1}{2}\right) -e_1 \right) \wp '\left(x+\frac{1}{2}\right).\nonumber
\end{align}
For the operators $A$ and $H$ there is an algebraic relation $A^2+Q(H)=0$
and the genus of this curve is $2$.

As in the previous example we calculate the coefficient of the monodromy (\ref{analcont}).
From the formulas $\wp(x)^2=\frac{1}{6}\wp''(x)+\frac{1}{12}g_2$, $\zeta '(x)=-\wp(x)$, $\zeta (x+1)-\zeta(x)= 2\eta_1$ and the periodicity of $\wp'(x)$, we have
\begin{align}
& \int_{0+\varepsilon }^{1+\varepsilon }\Xi (x,E)dx =E^2-(5e_1+8\eta_1)E+\frac{3}{4} g_2-27e_1^2+4e_1\eta_1 .\label{2100P1}
\end{align}
Since $E=-4e_1$ is a solution to $Q(E)=0$ such that $\Lambda (x+1,E_0)=\Lambda (x,E_0)$, it follows that
\begin{align} 
& \Lambda (x+1,E_*)=\Lambda (x,E_*) \exp \left( -\frac{1}{2} \int_{-4e_1}^{E_*}\frac{E^2-(5e_1+8\eta_1)E+\frac{3}{4} g_2-27e_1^2+4e_1\eta_1}{\sqrt{-Q(E)}} dE\right), \label{analcont2100} 
\end{align}
which is an hyperelliptic integral related to the genus two curve.

\section{Application to the spectral problem} \label{sec4}

\subsection{The strategy} $ $ \label{sec:str}
In this section eigenvalues of the Hamiltonian $H$ for the systems with physical boundary conditions are considered.
More precisely we consider eigenstates and eigenvalues of the Hamiltonian in an appropriate Hilbert space.
In the general setting the Hilbert space is taken as the space of functions that are square-integrable for some inner product.
For the infinite-dimensional case the Hamiltonian generally does not act on all elements of the whole Hilbert space but only on a dense subspace.
Hence we should specify the Hilbert space and its dense subspace for each model.

In section \ref{sec:Hilbsp} the Hilbert space and its dense subspace for the $BC_1$ Inozemtsev model are introduced. In previous sections we obtained the eigenfunctions $\Lambda(x,E)$ of the Hamiltonian $H$ in (\ref{eqn:Lam}) without boundary conditions and the constant for the quasi-periodicity in (\ref{analcont}).
It is mentioned in this section how the conditions for square-integrability correspond to the ones for quasi-periodicity.

For the trigonometric case, i.e. $\tau \rightarrow \sqrt{-1} \infty $, the model becomes much simpler. In this case all eigenvalues and eigenfunctions are described explicitly.
The eigenvalues of the nontrigonometric models are obtained by connecting from the trigonometric model.
Here the Hamiltonian $H$ depends upon the variable $p=\exp(\pi \sqrt{-1} \tau)$ and the trigonometric case corresponds to the case $p=0$.
It is shown in \cite{Tak2} that the eigenvalues are holomorphic in $p$ if $-1<p <1$ by an application of the theory of holomorphic perturbation due to Kato and Rellich \cite{Kat}.

To further understand the eigenvalues we examine the singular points of the eigenvalues in the variable $p$ on the complex domain $|p| <1$.
In section \ref{sec:condsq} a sufficient condition for the holomorphy of an eigenvalue at a given point $p=p_0$ is provided. Here relation (\ref{analcont}) plays an important role.

Firstly the case $l_0 \in \Zint_{\geq 1}$ and $l_1 \in \Zint_{\geq 1}$ is discussed in sections \ref{sec:Hilbsp} and \ref{sec:condsq}. 
In section \ref{sec:l10} we treat the case $l_0\in \Zint_{\geq 1}$ and $l_1 =0$, although the case $l_0=0$ and  $l_1 \in \Zint_{\geq 1}$ can be treated similarly.
In section \ref{sec:l0l10} the case $l_0=l_1=0$ is discussed.

\subsection{The Hilbert space and the trigonometric limit} $ $ \label{sec:Hilbsp}

In this subsection assume that $l_0 , l_1 \in \Zint _{\geq 1}$ and $l_2, l_3 \in \Zint _{\geq 0}$ and set $p=\exp(\pi \sqrt{-1} \tau)$. Then the potential of the Hamiltonian $H$ is dependent on $p$.
We introduce the Hilbert space and its dense subspace. Note that our notation differs from that in \cite{Tak2}. 
Define the following function, which is indeed the ground state for the trigonometric ($p=0$) model:
\begin{equation}
\Phi(x)=(\sin \pi x)^{l_0+1}(\cos \pi x)^{l_1+1}.
\label{Phix}
\end{equation}
By using the function $\Phi(x)$, we define the Hilbert space ${\bold H}$ and its subspace ${\bold D} $ for every $p$ ($|p|<1)$ as follows:
\begin{equation}
{\bold H} = \left\{ f\!  : \Rea \rightarrow \Cplx \mbox{ measurable} \left|
\begin{array}{l}
 \int_{0}^{1} |f(x)| ^2 dx<+\infty, \\
 f(x)=(-1)^{l_0+l_1} f(x+1) \mbox{ a.e. }x, \\
 f(x)=(-1)^{l_0+1} f(-x) \mbox{ a.e. }x
\end{array}
\right. \right\}
\label{Hilb1}
\end{equation}
and 
\begin{equation}
{\bold D} = \left\{ f\!  : \Rea \rightarrow \Cplx \left|
\begin{array}{l}
\frac{f(x)}{\Phi(x)} \mbox{ is holomorphic on }\Rea ,\\
 f(x)=(-1)^{l_0+l_1} f(x+1), \; f(x)=(-1)^{l_0+1} f(-x) 
\end{array}
\right. \right\}.
\label{Hilbdense1}
\end{equation}
Observe that the spaces ${\bold H}$ and ${\bold D}$ are independent of $p$.
It is shown that the space ${\bold D}$ is a dense subspace of $ {\bold H}$. 

\begin{remk}
Let 
$$
{\bold H}_{g} = \left\{ f\! : \Rea \rightarrow \Cplx  \mbox{ measurable} \left|
\begin{array}{l}
\int_{0}^{1} |f(x)| ^2 |\Phi (x)|^2 dx<+\infty, \\
 f(x)=f(x+1), \; f(x)=f(-x) \mbox{ a.e. }x  
\end{array}
\right. \right\}
$$
be the gauge-transformed Hilbert space. (This space is denoted as ${\bold H}$ in \cite{Tak2}.)
There is a simple relationship between the spaces ${\bold H}_{g}$ and ${\bold H}$, namely 
$
h(x) \in {\bold H}_{g} \; \Leftrightarrow h(x) \Phi(x) \in {\bold H}.
$
\end{remk}

It is shown in \cite{Tak2} that for each eigenstate $\psi _m(x) \in {\bold H}$ of the Hamiltonian $H$ we have $\psi _m(x) \in {\bold D}$ if $p\in \Rea $ and $|p|$ is sufficiently small. Thus it would be reasonable to consider the spectral problem for the space ${\bold D}$.

\begin{prop} \label{prop:HpresD}
The Hamiltonian $H$ preserves the space ${\bold D}$ for every $p\in \Cplx$ such that $|p|<1$.
\end{prop}
\begin{proof}
Let $p$ be a complex number such that $|p|<1$ and let ${\mathcal H}(p)=\Phi(x)^{-1} \circ H \circ \Phi(x) $ be the gauge-transformed Hamiltonian. Then we have 
\begin{equation}
{\mathcal H}(p)= -\frac{d^2}{dx^2}-2\pi\left( (l_0 +1) \frac{\cos \pi x}{\sin \pi x}- (l_1 +1)\frac{\sin \pi x}{\cos \pi x} \right) \frac{d}{dx}+V_p(x)
\end{equation}
 where $V_p(x)$ is an even periodic function in $x$ that is holomorphic on $x\in \Rea$ (for detais see \cite[\S 2]{Tak2}). It is enough to show that the operator ${\mathcal H}(p)$ preserves the space ${\bold D}_g$, where
\begin{equation}
{\bold D} _g= \left\{ f \! : \Rea \rightarrow \Cplx \left|
\begin{array}{l}
f(x) \mbox{ is holomorphic on }\Rea ,\\
 f(x)=f(x+1), f(x)= f(-x) 
\end{array}
\right. \right\}.
\label{Hilbdenseg1}
\end{equation}

Let $f(x) \in {\bold D}$. Then the function $\tilde{f}(x)={\mathcal H}(p)f(x)$ satisfies $\tilde{f}(x+1)=\tilde{f}(x)=\tilde{f}(-x)$ and $\tilde{f}(x-1/2)=\tilde{f}(1/2-x)$. The function $\tilde{f}(x)$ is meromorphic on $\Rea $ and may have poles at $x\in \frac{1}{2}\Zint $ of degree at most $1$.
If $\tilde{f}(x)$ has a pole at $x=0$ (resp. $x=1/2$), it contradicts the symmetry $\tilde{f}(x)=\tilde{f}(-x)$ (resp. $\tilde{f}(x-1/2)=\tilde{f}(1/2-x)$). Hence the function $\tilde{f}(x)$ has no poles at $x=0$ and $x=1/2$.
Furthermore from the periodicity the function $\tilde{f}(x)$ has no poles on $\Rea $ thus completing the proof. 
\end{proof}

For the case $p=0$ the problem becomes much easier. We briefly describe the eigenvalues (for details, see \cite{Tak2}).

If $p\rightarrow 0$ ($\tau \rightarrow \sqrt{-1}\infty$), then we have $\wp(x) \rightarrow -\pi^2/3+\pi ^2/\sin ^2\pi x$ etc. and $H \rightarrow H_T -\pi^2\left(\sum_{i=0}^3 l_i (l_i +1)\right) /3$, where 
\begin{align}
& H_T=  -\frac{d^2}{dx^2} + l_0(l_0+1)\frac{\pi ^2}{\sin ^2\pi x} +l_1(l_1+1)\frac{\pi ^2}{\cos ^2\pi x}. \label{triIno} 
\end{align}
The $(m+1)$st lowest eigenvalue of the trigonometric Hamiltonian $H_T$ on the space ${\bold H}$ is $\pi^2 (2m+l_0+l_1+2)^2$. The corresponding eigenfunction $\tilde{\psi} _m(x)$ is written as the product of the ground state $\Phi(x)$ from (\ref{Phix}) and the Jacobi polynomial of degree $m$. In particular we have $\tilde{\psi} _m(x) \in {\bold D}$.

If $-1<p<1$, then the existence, completeness, and real-holomorphy of the eigenvalues in $p$ are shown in \cite[Theorem 3.6]{Tak2} by application of Kato-Rellich theory.

\begin{prop} $($c.f. \cite[Theorem 3.6]{Tak2}$)$ \label{prop:complEm}
Let $l_0, l_1 \in \Zint_{\geq 1}$ and $l_2,l_3 \in \Zint_{\geq 0}$. All eigenvalues of $H$ for $-1<p<1$ on the Hilbert space ${\bold H}$ can be represented as $\tilde{E}_m \! (p)(=\pi^2 (2m+l_0+l_1+2)^2 -\frac{\pi^2}{3}\sum_{i=0}^3 l_i (l_i +1))$ $(m\in \Zint_{\geq 0})$, which is real-holomorphic in $p \in (-1,1)$, and the value $\tilde{E}_m\! (0)$ coincides with the $(m+1)$st lowest eigenvalue of the trigonometric model.
\end{prop}

Although the eigenvalue $\tilde{E}_m \! (p)$ is real-holomorphic, it may have singularities on the complex domain. This is investigated in the next subsection.

\subsection{The condition for square-integrability} $ $ \label{sec:condsq}

In this subsection we write down the condition to determine which values of $E$ admit square-integrable eigenfunctions and discuss the singularities of the eigenvalue $\tilde{E}_m\! (p)$ in the variable $p =\exp(\pi \sqrt{-1} \tau)$.
Throughout this section it is assumed that $l_0, l_1 \in \Zint_{\geq 1}$ and $l_2,l_3 \in \Zint_{\geq 0}$.

Let $\Lambda (x,E)$ be the function defined in (\ref{eqn:Lam}). Then $H\Lambda (x,E) = E \Lambda (x,E)$ and $\Lambda(x+1,E)= B(E) \Lambda (x,E)$, where $B(E)=\exp \left( -\frac{1}{2} \int_{E_0}^{E}\frac{Q_1(\tilde{E})}{\sqrt{-Q(\tilde{E})}} d\tilde{E}\right)$, $E_0$ is the value satisfying $\Lambda (x+1,E_0)=\Lambda (x,E_0)$ and $Q_1(E)= \int_{0+\varepsilon }^{1+\varepsilon }\Xi (x,E)dx$ (see Theorem \ref{thm:conj3}).

Set $\Lambda^{\mbox{\scriptsize sym}}(x,E)=\Lambda (x,E)-(-1)^{l_0}\Lambda(-x,E)$. Then the relation $H\Lambda^{\mbox{\scriptsize sym}} (x,E) = E\Lambda ^{\mbox{\scriptsize sym}} (x,E)$ is obvious. The following proposition gives a sufficient condition for the square-integrability because, if $\Lambda ^{\mbox{\scriptsize sym}} (x,E) \in {\bold D}$, then the function $\Lambda ^{\mbox{\scriptsize sym}} (x,E)$ is square-integrable on $(0,1)$.

\begin{prop} \label{prop:be1}
(i) If $B(E)=(-1)^{l_0+l_1}$, then $\Lambda ^{\mbox{\scriptsize sym}} (x,E) \in {\bold D}$. \\
(ii) If $Q(E) \neq 0$, then the function  $\Lambda ^{\mbox{\scriptsize sym}} (x,E)$ is nonzero.
\end{prop}
\begin{proof}
(i) We must check that the function $\Lambda ^{\mbox{\scriptsize sym}} (x,E)$ satisfies the defining relations (\ref{Hilbdense1}) for ${\bold D}$.

The relation $\Lambda ^{\mbox{\scriptsize sym}} (x,E) =(-1)^{l_0+1} \Lambda ^{\mbox{\scriptsize sym}} (-x,E)$ is trivial from the definition of $\Lambda ^{\mbox{\scriptsize sym}} (x,E)$.
The relation $\Lambda ^{\mbox{\scriptsize sym}} (x,E) =(-1)^{l_0+l_1} \Lambda ^{\mbox{\scriptsize sym}} (x+1,E)$ is obtained from the condition $B(E)=(-1)^{l_0+l_1}$.
Since the function  $\Lambda ^{\mbox{\scriptsize sym}} (x,E)$ satisfies the differential equation (\ref{InoEF}), it may have poles only at $x \in \frac{1}{2}\Zint $ on $\Rea$.

We have the expansion $\Lambda ^{\mbox{\scriptsize sym}} (x,E)= x^{\alpha }(c_0+c_1x+ \cdots )$, where $c_0 \neq 0$ and ($\alpha =-l_0$ or $l_0+1$), because the exponents of the differential equation (\ref{InoEF}) at $x=0$ are $-l_0$ and $l_0+1$. From the property $\Lambda ^{\mbox{\scriptsize sym}} (x,E) =(-1)^{l_0+1} \Lambda ^{\mbox{\scriptsize sym}} (-x,E)$ and $l_0 \in \Zint_{\geq 0}$, we have $\alpha =l_0+1$. Thus the function $\Lambda ^{\mbox{\scriptsize sym}} (x,E)$ is holomorphic at $x=0$.

We set $\tilde{\lambda}(x)=\Lambda ^{\mbox{\scriptsize sym}} (x-1/2,E)$. Then we have $ \tilde{\lambda}(x)= \Lambda ^{\mbox{\scriptsize sym}} (x-1/2,E) = (-1)^{l_0+l_1} \Lambda ^{\mbox{\scriptsize sym}} (x+1/2,E) = (-1)^{l_1+1} \Lambda ^{\mbox{\scriptsize sym}} (-x-1/2,E)= (-1)^{l_1+1} \tilde{\lambda}(-x)$. Similarly to the case $\Lambda ^{\mbox{\scriptsize sym}} (x,E)|_{x=0}$, the function $\tilde{\lambda}(x)$ is holomorphic at $x=0$. Hence the function $\Lambda ^{\mbox{\scriptsize sym}} (x,E)$ is holomorphic at $x=1/2$.

From the periodicity of $\Lambda ^{\mbox{\scriptsize sym}} (x,E)$ we obtain that the function  $\Lambda ^{\mbox{\scriptsize sym}} (x,E)$ is holomorphic on $\Rea$. 
Therefore we have checked all the defining relations and have proved (i).

(ii) It is sufficient to show that the functions $\Lambda (x,E)$ and $\Lambda(-x,E)$ are linearly independent. 
This is obvious from the expressions $\Lambda (-x,E)=\sqrt{\Xi (x,E)}\exp \int \frac{-\sqrt{-Q(E)}dx}{\Xi (x,E)}$ and $\Lambda (x,E)=\sqrt{\Xi (x,E)}\exp \int \frac{\sqrt{-Q(E)}dx}{\Xi (x,E)}$.
\end{proof}

\begin{prop} \label{prop:be2}
Let $l_0, l_1 \in \Zint_{\geq 1}$ and $l_2,l_3 \in \Zint_{\geq 0}$. If $B(E)\neq (-1)^{l_0+l_1}$ and  $Q(E)\neq 0$, then any nonzero solution to (\ref{InoEF}) is not square-integrable. In particular any value $E$ such that $B(E)\neq (-1)^{l_0+l_1}$ and  $Q(E)\neq 0$ is not an eigenvalue of the Hamiltonian $H$ in the space ${\bold D}$.
\end{prop}
\begin{proof}
From the condition $Q(E)\neq 0$ any solution to (\ref{InoEF}) can be written as a linear combination of $\Lambda (x,E)$ and $\Lambda (-x,E)$. From the inequality $B(E)\neq (-1)^{l_0+l_1}$ it is impossible to remove the poles at $x=0, \: 1/2$, and $1$ simultaneously. 
\end{proof}

For the case $Q(E)=0$ a different argument is required.

We observe the condition under which eigenvalues in ${\bold D}$ are analytically continued. 
The eigenvalues in ${\bold D}$ are determined by the condition $B(E)=(-1)^{l_0+l_1}$.
The function $B(E)$ depends on two variables, $E$ and $p$, and furthermore is holomorphic in $E$ and $p$ provided $(E,p) \in \Cplx \times B_1$ and $Q(E) \neq 0$, where $B_1=\{ p \in \Cplx \: | \: |p|<1 \}$.

Let $|p_{*}|<1$ and suppose that $B(E_*)=(-1)^{l_0+l_1}$ at $p=p_*$. By the implicit function theorem, if $\frac{\partial B}{\partial E} |_{p=p_* ,E=E_*} \neq 0$ and $Q(E_*) \neq 0$, then there exist $\varepsilon  \in \Rea_{>0}$ and a unique holomorphic function $E(p)$ on $\{ p\in \Cplx \: | \: |p- p_*|<\varepsilon \} $ such that $B(E(p))=(-1)^{l_0+l_1}$ and $E(p_*)=E_*$.
Thus we have the following proposition.
\begin{prop} \label{prop:analcontEp}
Let $|p_*|<1$ and $E_*$ be an eigenvalue of $H$ on ${\bold D}$ such that $\frac{\partial B}{\partial E} |_{p=p_* ,E=E_*} \neq 0$ and $Q(E_*) \neq 0$.
Then this eigenvalue on ${\bold D}$ is analytically continued near $p=p_*$ and it is holomorphic at $p=p_*$.
\end{prop}

From formula (\ref{analcont}) we have more explicit results.

\begin{thm} \label{thm:conti}
Let $l_0, l_1 \in \Zint_{\geq 1}$ and $l_2,l_3 \in \Zint_{\geq 0}$.\\
(i) If 
\begin{equation}
\int_{E_0}^{E_*}\frac{Q_1(E)}{\sqrt{-Q(E)}} dE \in 2(l_0+l_1) \pi\sqrt{-1} +4\pi \sqrt{-1} \Zint ,
\label{monodeqn}
\end{equation}
the value $E_0$ satisfies $\Lambda (x+1,E_0)= \Lambda (x,E_0)$, $Q_1(E)$ is a polynomial in $E$ defined by $Q_1(E)= \int_{0+\varepsilon }^{1+\varepsilon }\Xi (x,E)dx$ and $Q(E_*)\neq 0$, then we have $\Lambda ^{\mbox{\scriptsize sym}} (x,E_*) \in {\bold D}$ and $\Lambda ^{\mbox{\scriptsize sym}} (x,E_*) \neq 0$.\\
(ii) If $|p_*|<1$ and $E_*$ is such that $Q_1(E_*)\neq 0$, $Q(E_*)\neq 0$ and  there exists a nonzero function $f(x) \in {\bold D}$ satisfying $(H-E_*)f(x)=0$ at $p=p_*$, then the eigenvalue $E_*$ is analytically continued near $p=p_*$ and is holomorphic at $p=p_*$.
\end{thm}
\begin{proof}
From (\ref{analcont}) the condition $B(E_*)=(-1)^{l_0+l_1}$ is equivalent to the condition (\ref{monodeqn}).
In addition, as is shown in Propositions \ref{prop:be1} and \ref{prop:be2}, if $Q(E_*)\neq 0$, then the condition that there exists a nonzero function $f(x) \in {\bold D}$ such that $(H-E_*)f(x)=0$  is also equivalent to $B(E_*) =(-1)^{l_0+l_1}$. 

Finally it is easy to see that the conditions $\frac{\partial B}{\partial E} |_{p=p_* ,E=E_*} \neq 0$ and $Q(E_*) \neq 0$ are equivalent to the conditions  $Q_1(E_*)\neq 0$ and $Q(E_*)\neq 0$.
By combining these facts and Proposition \ref{prop:analcontEp}, we obtain the theorem.
\end{proof}
\begin{remk}
If the value $E_0$ satisfies $\Lambda (x+1,E_0)= -\Lambda (x,E_0)$, then the condition (\ref{monodeqn}) is replaced by 
\begin{equation}
\int_{E_0}^{E_*}\frac{Q_1(E)}{\sqrt{-Q(E)}} dE \in 2(l_0+l_1+1) \pi \sqrt{-1} +4\pi  \sqrt{-1} \Zint.
\end{equation}
\end{remk}

Let $\tilde{E}_m \! (p)$ be the eigenvalue of the Hamiltonian $H$ on the Hilbert space ${\bold H}$ defined in Proposition \ref{prop:complEm}.
Then from \cite{Tak2} the eigenfunction $v_m(x,p)$ of the eigenvalue $\tilde{E}_m \! (p)$ on the full Hilbert space ${\bold H}$ satisfies $v_m(x,p) \in {\bold D}$, if $p \in \Rea $ and $|p|$ is sufficiently small.
Hence $\int_{E_0}^{\tilde{E}_m \! (p)}\frac{Q_1(E)}{\sqrt{-Q(E)}} dE \in 2(l_0+l_1) \pi \sqrt{-1} +4\pi \sqrt{-1} \Zint$, if $E_0$ satisfies $\Lambda (x+1,E_0)= \Lambda (x,E_0)$, $p \in \Rea $ and $|p|$ is sufficiently small. From Proposition \ref{prop:complEm} the eigenvalue $\tilde{E}_m \! (p)$ is holomorphic if $-1<p<1$.

From Theorem \ref{thm:conti}, if $Q(\tilde{E}_m \! (p))\neq 0$ and $E_0$ satisfies $\Lambda (x+1,E_0)= \Lambda (x,E_0)$, then the eigenvalue $\tilde{E}_m \! (p)$ is so continued analytically in $p \: (\in \Cplx $) as to maintain the relation $\int_{E_0}^{\tilde{E}_m \! (p)}\frac{Q_1(E)}{\sqrt{-Q(E)}} dE \in 2(l_0+l_1) \pi \sqrt{-1} +4\pi \sqrt{-1} \Zint$. Let $p_* \in \Cplx$ such that $|p_*|<1$.
If $Q(\tilde{E}_m \! (p_*))\neq 0$ and $Q_1(\tilde{E}_m \! (p_*))\neq 0$, then the eigenvalue $\tilde{E}_m \! (p)$ does not have a singularity at $p=p_*$;
in other words the eigenvalue $\tilde{E}_m \! (p)$ can have a singularity at $p=p_*$ only when $Q(\tilde{E}_m \! (p_*))= 0$ or $Q_1(\tilde{E}_m \! (p_*))= 0$.

As an aside note that in \cite{Tak2} an algorithm for calculating the eigenvalue $\tilde{E}_m \! (p)$ as a power series in $p$ was introduced and the convergence of the power series for sufficiently small $|p|$ was shown.
Numerically it can be demonstrated that the radius of convergence of the power series is strictly less than $1$. Assume that the radius of convergence of the power series is $p_{c} \in (0,1)$. 
Then from the preceding discussion there exists a solution to $Q(\tilde{E}_m \! (p_*))= 0$ or $Q_1(\tilde{E}_m \! (p_*))= 0$ satisfying $|p_*|=p_{c}$.
In this way the condition $(Q(E)= 0$ or $Q_1(E)= 0)$ is closely connected to the radius of convergence of the eigenvalues.

For the case $l_0=2$, $l_1=1$, $l_2=0$ and $l_3=0$ the polynomial $Q(E)$ is given in (\ref{2100P}) and the polynomial $Q_1(E)$ is given as the right-hand side of (\ref{2100P1}).
Then we have $\deg _E Q(E)=5$ and $\deg _E Q_1(E)=2$.
To find the singular points of the eigenvalues for the space ${\bold D}$ it would be effective to compute the values $E_*$ such that ($Q(E_*)=0$ or $Q_1(E_*)=0$) as $p_*$ varies and investigate when these values $E_*$ also satisfy $\int_{-4e_1}^{E_*}\frac{Q_1(E)}{\sqrt{-Q(E)}} dE \in (4\Zint +2)\pi \sqrt{-1} $.
It is anticipated that more such results will be reportable in the future.

\subsection{The case $l_0 \in \Zint_{\geq 1}$ and $l_1=0$} $ $ \label{sec:l10}

In this subsection we consider the case $l_0 \in \Zint_{\geq 1}$, $l_1=0$, $l_2 ,l_3 \in \Zint_{\geq 0}$.
For this case consider the larger Hilbert space ${\bold H}^{\star}$ and its dense subspace ${\bold D}^{\star}$ defined by
\begin{equation}
{\bold H}^{\star} = \left\{ f \! : \Rea \rightarrow \Cplx \mbox{ measurable} \left|
\begin{array}{l}
 \int_{0}^{1} |f(x)| ^2 dx<+\infty, \\
 f(x)= f(x+2) \mbox{ a.e. }x, \\
 f(x)=(-1)^{l_0+1} f(-x) \mbox{ a.e. }x
\end{array}
\right. \right\}
\label{Hilb1l10}
\end{equation}
and 
\begin{equation}
{\bold D}^{\star} = \left\{ f\! : \Rea \rightarrow \Cplx \left|
\begin{array}{l}
\frac{f(x)}{(\sin \pi x )^{l_0 +1}} \mbox{ is holomorphic on }\Rea ,\\
 f(x)= f(x+2), \; f(x)=(-1)^{l_0+1} f(-x)
\end{array}
\right. \right\}.
\label{Hilbdense1l10}
\end{equation}
Then ${\bold D}^{\star}\subset {\bold H}^{\star}$, ${\bold H} \subset {\bold H}^{\star}$, and ${\bold D} \subset {\bold D}^{\star}$.
It can be shown that the Hamiltonian $H$ preserves the space ${\bold D}^{\star}$ (c.f. Proposition \ref{prop:HpresD}).

For the case $l_0 \in \Zint_{\geq 1}$, $l_1=0$, $l_2 ,l_3 \in \Zint_{\geq 0}$, all eigenvalues of $H$ for $-1<p<1$ on the Hilbert space ${\bold H}^{\star}$ can be represented as $\tilde{E}_m ^{\star} (p)$ $(m\in \Zint_{\geq 0})$, which is real-holomorphic in $p \in (-1,1)$ and for which $\tilde{E}_m^{\star} (0)( =\pi ^2 (m+l_0+1)^2 -\frac{\pi ^2}{3}\sum_{i=0}^3 l_i (l_i +1))$ coincides with the $(m+1)$st lowest eigenvalue of the trigonometric model (see \cite[Theorem 3.6]{Tak2}).

Let $\Lambda (x,E)$ be the function defined in (\ref{eqn:Lam}). Then $H\Lambda (x,E) = E\Lambda (x,E)$ and $\Lambda(x+1,E)= B(E) \Lambda (x,E)$, where $B(E)= \exp \left( -\frac{1}{2} \int_{E_0}^{E}\frac{Q_1(\tilde{E})}{\sqrt{-Q(\tilde{E})}} d\tilde{E} \right)$, $E_0$ satisfies $\Lambda (x+1,E_0)= \Lambda (x,E_0)$, and $Q_1(E)= \int_{0+\varepsilon }^{1+\varepsilon }\Xi (x,E)dx$.
Set $\Lambda^{\mbox{\scriptsize sym}}(x,E)=\Lambda (x,E)-(-1)^{l_0}\Lambda(-x,E)$. 
If $B(E)=1$ or $-1$, then we have $\Lambda ^{\mbox{\scriptsize sym}} (x,E) \in {\bold D}^{\star}$ and, if $Q(E) \neq 0$, then the function  $\Lambda ^{\mbox{\scriptsize sym}} (x,E)$ is nonzero.
Conversely, if $B(E)\neq \pm 1$ and $Q(E)\neq 0$, then it can be shown that any nonzero solution to (\ref{InoEF}) is not square-integrable.
Note that the condition $B(E)= \pm 1$ is equivalent to the condition $\int_{E_0}^{E}\frac{Q_1(\tilde{E})}{\sqrt{-Q(\tilde{E})}} d\tilde{E} \in 2\pi \sqrt{-1} \Zint$, where $E_0$ satisfies $Q(E_0)=0$.

Similar statements to Proposition \ref{prop:analcontEp} and Theorem \ref{thm:conti} can also be derived.

\subsection{The case $l_0 = l_1=0$} $ $ \label{sec:l0l10}

In this subsection we consider the case $l_0 =l_1=0$, $l_2 ,l_3 \in \Zint_{\geq 0}$.

Set
\begin{equation}
{\bold H}^{\star \star} = \left\{ f\! : \Rea \rightarrow \Cplx \mbox{ measurable} \left|
\begin{array}{l}
 \int_{0}^{1} |f(x)| ^2 dx<+\infty, \\
 f(x)= f(x+2) \mbox{ a.e. }x
\end{array}
\right. \right\}
\label{Hilb100}
\end{equation}
and 
\begin{equation}
{\bold D}^{\star \star} = \left\{ f\! : \Rea \rightarrow \Cplx \left|
\begin{array}{l}
f(x) \mbox{ is holomorphic on }\Rea ,\\
 f(x)= f(x+2)
\end{array}
\right. \right\}.
\label{Hilbdense00}
\end{equation}
Then ${\bold D}^{\star \star }  \subset{\bold H}^{\star \star }$, ${\bold H} \subset {\bold H}^{\star } \subset {\bold H}^{\star \star }$, and ${\bold D} \subset {\bold D}^{\star} \subset {\bold D}^{\star \star }$.

It is shown that the Hamiltonian $H$ preserves the space ${\bold D}^{\star \star }$ (c.f. Proposition \ref{prop:HpresD}).

For the case $l_0 =l_1=0$, $l_2 ,l_3 \in \Zint_{\geq 0}$, all eigenvalues of $H$ for $-1<p<1$ on the Hilbert space ${\bold H}^{\star \star }$ can be represented as $\tilde{E}_m ^{\diamondsuit} (p)$ $(m\in \Zint_{\geq 0})$ and $\tilde{E}_m^{\Box} (p)$  $(m\in \Zint_{\geq 1})$, which are real-holomorphic in $p \in (-1,1)$ and for which $\tilde{E}_m^{\diamondsuit} (0) =\pi^2 m^2 -\frac{\pi^2}{3}\sum_{i=2}^3 l_i (l_i +1)$ $(m\in \Zint_{\geq 0})$ and $\tilde{E}_m^{\Box} (0) =\pi^2 m^2 -\frac{\pi^2}{3}\sum_{i=2}^3 l_i (l_i +1)$ $(m\in \Zint_{\geq 1})$ (see \cite[Theorem 3.6]{Tak2}).

Let $\Lambda (x,E)$ be the function defined in (\ref{eqn:Lam}). Then $H\Lambda (x,E) = E\Lambda (x,E)$ and $\Lambda(x+1,E)= B(E) \Lambda (x,E)$, where $B(E)= \exp\left( -\frac{1}{2}\int_{E_0}^{E}\frac{Q_1(E)}{\sqrt{-Q(E)}} dE \right)$, $E_0$ satisfies $\Lambda (x+1,E_0)= \Lambda (x,E_0)$ and $Q_1(E)= \int_{0+\varepsilon }^{1+\varepsilon }\Xi (x,E)dx$.
It is easy to show that, if $B(E)=\pm 1$, then $\Lambda (x,E) \in {\bold D}^{\star \star }$.
The condition  $B(E_*)=\pm 1$ is equivalent to the condition $\int_{E_0}^{E_*}\frac{Q_1(E)}{\sqrt{-Q(E)}} dE \in 2\pi \sqrt{-1} \Zint$.

For the case $l_0 =l_1=0$, $l_2 ,l_3 \in \Zint_{\geq 0}$ we have the following proposition.
Some parts are proved similarly to Proposition \ref{prop:analcontEp} and other parts are proved directly.

\begin{prop} \label{prop:00}
(i) If $Q(E_*)\neq 0$ and $B(E_*)=\pm 1$, then the functions $\Lambda (x,E_*)$ and  $\Lambda (-x,E_*)$ are linearly independent and thus $E_*$ is a double eigenvalue.\\
(ii) If $Q(E_*)= 0$, then the function $\Lambda (x,E_*)$ is doubly periodic up to signs and the relations $(\Lambda (x,E_*)=\Lambda (-x,E_*)$ or $-\Lambda (-x,E_*))$ and $B(E_*)=\pm 1$ are obtained. In this case the eigenfunction $\Lambda (x,E_*)$ corresponds to the nonrepeated eigenvalue $E_*$.\\
(iii) The nonrepeated eigenvalues are obtained by the diagonalization of the finite-dimensional space $V$ defined in section \ref{sec:invsp}.\\ 
(iv) The repeated eigenvalues $E_*$ are analytically so continued in $p$ as to preserve the property $\int_{E_0}^{E_*}\frac{Q_1(E)}{\sqrt{-Q(E)}} dE \in 2\pi \sqrt{-1} \Zint$, where $E_0$ satisfies $Q(E_0)=0$. If $Q_1(E_*)\neq 0$ and $Q(E_*)\neq 0$ at $p=p_*$, then the eigenvalue is analytic in $p$ at $p=p_*$.
\end{prop}

By application of Proposition \ref{prop:00} the following statements related to the finite-gap property can be shown:\vspace{3mm}
\\
${\bullet}$ If $l_2+l_3$ is even and $-1<p<1$, then $\tilde{E}_{2m}^{\diamondsuit} (p)=\tilde{E}_{2m}^{\Box} (p)$ for $m > \frac{1}{2}(l_2+l_3)$ and $\tilde{E}_{2m-1}^{\diamondsuit} (p)=\tilde{E}_{2m-1}^{\Box} (p)$ for $m > \frac{1}{2}|l_2-l_3|$. The rest of the eigenvalues correspond to the ones of the Hamiltonian $H$ on the space $V$ defined in section \ref{sec:invsp}.\\
${\bullet}$ If $l_2+l_3$ is odd and $-1<p<1$, then $\tilde{E}_{2m}^{\diamondsuit} (p)=\tilde{E}_{2m}^{\Box} (p)$ for $m > \frac{1}{2}(|l_2-l_3|-1)$ and $\tilde{E}_{2m-1}^{\diamondsuit} (p)=\tilde{E}_{2m-1}^{\Box} (p)$ for $m > \frac{1}{2}(l_2+l_3+1)$. The rest of the eigenvalues correspond to the ones of the Hamiltonian $H$ on the space $V$ defined in section \ref{sec:invsp}.\\

For the case $l_0=l_1=l_2=0$ and $l_3=1$ the results are illustrated more explicitly.
In this case $Q(E)=(E+e_1)(E+e_2)(E+e_3)$ and $Q_1(E)=E-2\eta_1$. 
From Proposition \ref{prop:00} the nonrepeated eigenvalues are $-e_1, -e_2$ and $-e_3$. Moreover $\tilde{E}_0 ^{\diamondsuit} (p)=-e_1$ and $\{ \tilde{E}_1^{\Box} (p), \: \tilde{E}_1 ^{\diamondsuit} (p)\} =\{ -e_2, -e_3\}$ whereas, for the other eigenvalues, $\tilde{E}_m^{\Box} (p)= \tilde{E}_m ^{\diamondsuit} (p)$ ($m \in \Zint_{\geq 2}$).
Finally the double eigenvalue $\tilde{E}_m^{\Box} (p)$ ($m \in \Zint_{\geq 2}$) is analytically so continued as to preserve $\int_{-e_1}^{\tilde{E}_m^{\Box} (p)}\frac{E-2\eta_1}{\sqrt{-(E+e_1)(E+e_2)(E+e_3)}} dE \in 2\pi \sqrt{-1} \Zint$ and, if $\tilde{E}_m^{\Box} (p_*) \neq 2\eta_1, -e_1, -e_2$ or $ -e_3$, then $\tilde{E}_m^{\Box} (p)$ is analytic at $p=p_*$.

\section{The Heun equation and the Inozemtsev model} \label{sec:Heun}

In this subsection the relationship between the Heun equation and the Inozemtsev model is reviewed and the results in this paper are translated to the ones for the Heun equation.

Let $f(x)$ be an eigenfunction of the Hamiltonian $H$ with the eigenvalue $E$ (\ref{InoEF}), i.e.,  
\begin{equation} \left( -\frac{d^2}{dx^2} + \sum_{i=0}^3 l_i(l_i+1)\wp (x+\omega_i)-E\right) f(x)=0.
\label{InoEF2}
\end{equation}
Set
\begin{equation}
a=\frac{e_2-e_3}{e_1-e_3}, \; \; \; \; \tilde{\Phi}(w)=w^{\frac{l_0+1}{2}}(w-1)^{\frac{l_1+1}{2}}(aw-1)^{\frac{l_2+1}{2}}.
\label{tilP}
\end{equation}
Under the transformation 
\begin{equation}
w=\frac{e_1-e_3}{\wp (x)-e_3}
\label{wxtrans}
\end{equation}
equation (\ref{InoEF2}) is equivalent to the equation
\begin{equation}
\left( \! \left(\frac{d}{dw}\right) ^2 \! + \left( \frac{\gamma}{w}+\frac{\delta }{w-1}+\frac{\epsilon}{w-t}\right) \frac{d}{dw} +\frac{\alpha \beta w -q}{w(w-1)(w-t)} \right)\tilde{f}(w)=0, 
\label{Heun}
\end{equation}
where $\tilde{f}(\frac{e_1-e_3}{\wp (x)-e_3})\tilde{\Phi}(\frac{e_1-e_3}{\wp (x)-e_3})=f(x)$, $t=\frac{1}{a}$, $\alpha =\frac{1}{2}(l_0+l_1+l_2+l_3+4)$, $\beta =\frac{1}{2}(l_0+l_1+l_2-l_3+3)$, $\gamma = l_0+\frac{3}{2}$, $\delta = l_1+\frac{3}{2}$, $\epsilon =l_2+\frac{3}{2}$, $q=-\frac{1}{4a}\left( \frac{E}{e_1-e_3}+c_0 \right)$, and $c_0= \frac{a+1}{3}\sum_{i=0}^3 l_i(l_i +1)-a(l_0+l_2+2)^2-(l_0+l_1+2)^2$ (for details see \cite{Tak2}). Note that the relation $\gamma +\delta +\epsilon =\alpha +\beta +1$ is satisfied. 

Equation (\ref{Heun}) with the condition $\gamma +\delta +\epsilon =\alpha +\beta +1$ is called the Heun equation. The Heun equation is the standard canonical form of a Fuchsian equation with four singularities.
By the aforementioned correspondence between $x$ and $w$, the parameters $(l_0, l_1, l_2, l_3, E, \tau)$ correspond essentially one-to-one to the parameters $( \alpha ,\beta ,\gamma ,\delta ,\epsilon, q, t)$ with the condition that $\gamma +\delta +\epsilon =\alpha +\beta +1$ (see \cite{Ron,OOS,Tak2}).

\begin{remk}
The expression of the Heun equation in terms of elliptic functions as in (\ref{InoEF2}) was essentially established by Darboux more than 100 years ago (see \cite{Smi}).

There are other possible choices for the function $\tilde{\Phi}(w)$ of (\ref{tilP}) and the relation (\ref{wxtrans}) between $x$ and $w$.
For example, Smirnov \cite{Smi} chose the correspondence defined by $t= \frac{e_3-e_1}{e_2-e_1}$, $\tilde{\Phi}(w)=w^{-l_1/2}(w-1)^{-l_2/2}(w-t)^{-l_3/2}$, and $w=\frac{\wp(x) -e_1}{e_2-e_1}$.
\end{remk}

The transformation $w=\frac{e_1-e_3}{\wp (x)-e_3}$ produces the following correspondence:
\begin{equation}
\begin{array}{|c||c|c|c|c|}
\hline
 x & 0 & \frac{1}{2} & \frac{\tau+1}{2} & \frac{\tau}{2} \\ 
\hline
 w & 0 & 1 & \frac{1}{a} & \infty \\ 
\hline
\end{array}
\end{equation}
Other mathematical objects and concepts discussed in this paper are transformed under (\ref{wxtrans}) as follows:\vspace{3mm}
\\
$\bullet$ In sections \ref{sec2} and \ref{sec3}, we considered the finite-gap property for the Inozemtsev model for the case $l_0 ,l_1, l_2, l_3 \in \Zint $.
The condition $l_0 ,l_1, l_2, l_3 \in \Zint $ is transformed into the condition $\gamma, \delta, \epsilon , \alpha-\beta \in \frac{1}{2}+\Zint $.\\
$\bullet$ Some statements in this paper have assumptions such as $l_i =0$. 
If the assumption $l_i=0$ is changed to $l_i=-1$, all the applicable statements in this paper remain valid.
The condition $l_0=0$ or $-1$ (resp. $l_1=0$ or $-1$, $l_2=0$ or $-1$, $l_3=0$ or $-1$) is transformed under (\ref{wxtrans}) into the condition $\gamma =\frac{1}{2}$ or $\frac{3}{2}$ (resp. $\delta =\frac{1}{2}$ or $\frac{3}{2}$, $\epsilon =\frac{1}{2}$ or $\frac{3}{2}$, $\alpha -\beta= \pm \frac{1}{2}$).\\
$\bullet$ In section \ref{sec:invsp} we introduced the doubly periodic functions up to signs. The condition of doubly periodicity up to signs in $x$ is transformed into algebraicity in $w$.\\
$\bullet$ By translating the functions (e.g., $\Lambda (x,E)$, $\Xi(x,E)$) and the formulas (e.g., (\ref{const}), (\ref{eqn:Lam}) and (\ref{eqn:Lam01})) in this paper or in \cite{Tak1} into a form appropriate for the Heun equation (\ref{Heun}), we recover the corresponding formulas of Smirnov \cite{Smi}.\\
$\bullet$ Let ${\mathcal C}$ be the cycle enclosing two points $0$ and $1$ in $w$.
The path from $\varepsilon $ to $\varepsilon +1$ ($|\varepsilon |$ sufficiently small) in $x$ corresponds to the cycle ${\mathcal C}$ in $w$ (see Figure 1).\vspace{3mm}
\begin{figure}
\begin{picture}(380,113)(10,0)
\put(10,20){\circle*{2}}
\put(100,20){\circle*{2}}
\put(10,90){\circle*{2}}
\put(13,12){$0$}
\put(13,82){$\tau$}
\put(103,12){$1$}
\put(18,110){$x$}
\put(15,105){\line(1,0){13}}
\put(28,105){\line(0,1){13}}
\put(10,35){\vector(1,0){90}}
\put(150,5){\line(1,6){19}}
\put(190,40){\circle*{2}}
\put(260,40){\circle*{2}}
\put(350,40){\circle*{2}}
\put(193,32){$0$}
\put(353,32){$\frac{1}{a}$}
\put(263,32){$1$}
\put(243,61){${\mathcal C}$}
\put(178,110){$w$}
\put(175,105){\line(1,0){13}}
\put(188,105){\line(0,1){13}}
\put(195,55){\line(1,0){65}}
\put(260,40){\oval(30,30)[r]}
\put(260,25){\line(-1,0){70}}
\put(190,40){\oval(30,30)[l]}
\put(190,55){\vector(1,0){3}}
\end{picture}
\caption{}
\end{figure}

Assume $\gamma, \delta, \epsilon , \alpha-\beta \in \frac{1}{2}+\Zint $ and $\gamma +\delta +\epsilon =\alpha +\beta +1$.
The coefficients of the monodromy on the cycle ${\mathcal C}$ are written as the hyperelliptic integral in $q$ (the parameter defined after equation (\ref{Heun})).
Concretely let $\tilde{\Lambda}_1(w,q)$ and $\tilde{\Lambda}_2(w,q)$ be the functions defined by the relations $\tilde{\Lambda}_1\left(\frac{e_1-e_3}{\wp(x) -e_3}, -\frac{1}{4a}\left(\frac{E}{e_1-e_3}+c_0\right)\right) \tilde{\Phi}\left( \frac{e_1-e_3}{\wp(x) -e_3}\right) = \Lambda (x,E)$ and $\tilde{\Lambda}_2\left(\frac{e_1-e_3}{\wp(x) -e_3}, -\frac{1}{4a}\left(\frac{E}{e_1-e_3}+c_0\right)\right) \tilde{\Phi}\left( \frac{e_1-e_3}{\wp(x) -e_3}\right) = \Lambda (-x,E)$. Then the functions $\tilde{\Lambda}_1(w,q)$ and $\tilde{\Lambda}_2(w,q)$ satisfy equation (\ref{Heun}).
Let $\tilde{Q} (q)$ and $\tilde{Q}_1(q)$ be the monic polynomials defined by $Q(E)=(4(e_3-e_2))^{2g+1}\tilde{Q}\left(\frac{-1}{4a}\left( \frac{E}{e_1-e_3}+c_0\right)\right)$ and $Q_1(E)=(4(e_3-e_2))^{g}\tilde{Q}_1\left( \frac{-1}{4a}\left(\frac{E}{e_1-e_3}+c_0\right)\right)$,
fix $q_*$ and assume that $\tilde{Q}(q_*)\neq 0$. Finally let $\tilde{\Lambda}^{{\mathcal C}}_1(w,q_*)$ and $\tilde{\Lambda}^{{\mathcal C}}_2(w,q_*)$ be the functions obtained by analytic continuation in $w$ along the cycle ${\mathcal C}$ from the functions $\tilde{\Lambda}_1(w,q_*)$ and $\tilde{\Lambda}_2(w,q_*)$, and define the monodromy matrix $M^{{\mathcal C}}$ by 
\begin{equation}
(\tilde{\Lambda}^{{\mathcal C}}_1(w,q_*), \tilde{\Lambda}^{{\mathcal C}}_2(w,q_*))=(\tilde{\Lambda}_1(w,q_*), \tilde{\Lambda}_2(w,q_*))M^{{\mathcal C}} .
\end{equation}
Let $q_0$ be a solution to the equation $\tilde{Q}(q_0)=0$. Then we have $\tilde{\Lambda}_1(w,q_0) =\tilde{\Lambda}_2(w,q_0)$ and $\tilde{\Lambda}_1^{{\mathcal C}} (w,q_0) = \pm \tilde{\Lambda}_1(w,q_0)$ (c.f. Corollary \ref{cor:conj3}).
From formula (\ref{analcont}) $M^{{\mathcal C}}$ is a diagonal matrix and the eigenvalues of $M^{{\mathcal C}}$ are given by the following functions:
\begin{align}
& \pm (-1)^{l_0+l_1} \exp \left( \sqrt{(e_3-e_2)} \int_{q_0}^{q_*}\frac{\tilde{Q}_1(q)}{\sqrt{-\tilde{Q}(q)}} dq\right),
\label{analcontq} \\
& \pm (-1)^{l_0+l_1} \exp \left( - \sqrt{(e_3-e_2)} \int_{q_0}^{q_*}\frac{\tilde{Q}_1(q)}{\sqrt{-\tilde{Q}(q)}} dq\right), \nonumber
\end{align}
where the sign $\pm$ is chosen to be the same as the one $\tilde{\Lambda}_1^{{\mathcal C}} (w,q_0) = \pm \tilde{\Lambda}_1(w,q_0)$.
Note that the factor $(-1)^{l_0+l_1}$ comes from the branch of the function $\tilde{\Phi}(w)$.

For the case $\tilde{Q}(q_*)=0$ it is necessary to select the other basis for the solutions to (\ref{Heun}) which entails other calculations.

\appendix
\section {Elliptic functions} \label{sect:append}
This appendix presents the definitions of and formulas for elliptic functions.

Let $\omega_1$ and $\omega_3$ be complex numbers such that the value $\omega_3/ \omega_1$ is an element of the upper half plane. In this paper we set $\omega_1=1/2$ and $\omega_3=\tau/2$.

The Weierstrass $\wp$-function, the Weierstrass sigma-function and the Weierstrass zeta-function with periods $(2\omega_1, 2\omega_3)$ are defined as follows:

\begin{align}
& \wp (z)= \frac{1}{z^2}+  \sum_{(m,n)\in \Zint \times \Zint \setminus \{ (0,0)\} } \left( \frac{1}{(z-2m\omega_1 -2n\omega_3)^2}-\frac{1}{(2m\omega_1 +2n\omega_3)^2}\right),  \\
& \sigma (z)=z\prod_{(m,n)\in \Zint \times \Zint \setminus \{(0,0)\} } \left(1-\frac{z}{2m\omega_1 +2n\omega_3}\right) \nonumber \\
& \; \; \; \; \; \; \; \; \; \; \; \; \; \; \cdot \exp\left(\frac{z}{2m\omega_1 +2n\omega_3}+\frac{z^2}{2(2m\omega_1 +2n\omega_3)^2}\right), \nonumber \\
& \zeta(z)=\frac{\sigma'(z)}{\sigma (z)}. \nonumber
\end{align}
Setting $\omega_2=-\omega_1-\omega_3$ and 
\begin{align}
& e_i=\wp(\omega_i), \; \; \; \eta_i=\zeta(\omega_i) \; \; \; \; (i=1,2,3)
\end{align}
yields the relations
\begin{align}
& e_1+e_2+e_3=\eta_1+\eta_2+\eta_3=0, \; \; \; \\
& \wp(z)=-\zeta'(x), \; \; \; (\wp'(z))^2=4(\wp(z)-e_1)(\wp(z)-e_2)(\wp(z)-e_3), \nonumber \\
& \wp(z+2\omega_j)=\wp(z), \; \; \; \zeta(z+2\omega_j)=\zeta(z)+2\eta_j \; \; \; \; (j=1,2,3), \nonumber \\
& \wp(z+\omega_i)=e_i+\frac{(e_i-e_{i'})(e_i-e_{i''})}{\wp(z)-e_i} \; \; \; \; (i=1,2,3),\nonumber
\end{align}
where $i', i'' \in \{1,2,3\}$ with $i'<i''$, $i\neq i'$, and $i\neq i''$.
The constants $g_2$ and $g_3$ are defined by
\begin{equation}
g_2=-4(e_1e_2+e_2e_3+e_3e_1), \; \; \; g_3=4e_1e_2e_3.
\end{equation}
The co-sigma functions $\sigma_i(z)$ $(i=1,2,3)$ and co-$\wp$ functions $\wp_i(z)$ $(i=1,2,3)$ are defined by
\begin{align}
& \sigma_i(z)=\exp (-\eta_i z)\frac{\sigma(z+\omega_i)}{\sigma(\omega _i)}, \; \; \; \wp_i(z) = \frac{\sigma_i(z)}{\sigma(z)},
\end{align}
and satisfy
\begin{align}
& \wp_i(z) ^2 =\wp(z)-e_i \; \; \; \; \; \; \; (i=1,2,3). \label{sigmai}
\end{align}

{\bf Acknowledgment}
The author thanks Professors ~O. Chalykh, ~E. Date, ~Y. Komori , ~V. Kuznetsov, ~K. Takasaki, ~T. Takebe and ~A. Veselov for fruitful discussions. 
He was partially supported by the Grant-in-Aid for Scientific Research (No. 13740021) from the Japan Society for the Promotion of Science.

\label{lastpage}

\end{document}